\documentclass[a4paper,11pt,fleqn]{article}
\usepackage{amsmath,amssymb,amsthm,graphicx,subfigure,float,caption,epstopdf,tabularx,color, bm,amsfonts,epic}
\usepackage[top=1in, bottom=1in, left=1.25in, right=1.25in]{geometry}
\usepackage{appendix}
\usepackage{multirow}
\usepackage{diagbox}
\usepackage{appendix}
\allowdisplaybreaks
\newtheorem{thm}{Theorem}[section]

\newtheorem{rmk}{Remark}[section]
\newtheorem{PPS}{Proposition}[section]
\newtheorem{shm}{Scheme}[section]

\newtheorem*{prf}{Proof}
\numberwithin{equation}{section}

\usepackage{indentfirst}
\graphicspath{{Fig}}

\begin{document}
\title{High-order linearly implicit structure-preserving exponential integrators for the nonlinear Schr\"odinger equation}
\author{Chaolong Jiang$^{1,2}$, Jin Cui$^{3}$, Xu Qian$^1$\footnote{Correspondence author. Email:
qianxu@nudt.edu.cn.},  and Songhe Song$^1$ \\
{\small $^1$ Department of Mathematics, College of Liberal Arts and Science, }\\
{\small National University of Defense Technology, Changsha, 410073, P.R. China}\\
{\small $^2$ School of Statistics and Mathematics, }\\
{\small Yunnan University of Finance and Economics, Kunming 650221, P.R. China}\\
{\small $^3$ Department of Basic Sciences, Nanjing Vocational College of Information Technology,}\\
  {\small Nanjing 210023, China} \\
}
\date{}
\maketitle

\begin{abstract}
A novel class of high-order linearly implicit energy-preserving  integrating factor Runge-Kutta methods are proposed for the nonlinear Schr\"odinger equation. { Based on the idea of the scalar auxiliary variable approach, the original equation is first reformulated into an equivalent form which satisfies a quadratic energy}. The spatial derivatives of the system are
then approximated with the standard Fourier pseudo-spectral method. Subsequently, we apply the extrapolation technique/prediction-
correction strategy to the nonlinear terms of the semi-discretized system and a linearized energy-conserving system is obtained. A fully discrete scheme is gained by further using the integrating factor Runge-Kutta method to the resulting system. We show that, under certain circumstances for the coefficients of a Runge-Kutta method, the proposed scheme can produce numerical solutions along which the modified energy is precisely conserved, as is the case with the analytical solution and is extremely efficient in the sense that only linear equations with constant coefficients need to be solved at every time step. Numerical results are addressed to demonstrate the remarkable superiority of the proposed schemes in comparison with other existing structure-preserving schemes.  \\[2ex]
\end{abstract}

\section{Introduction}
We consider the nonlinear Schr\"odinger equation (NLSE) given as follows:
\begin{align}\label{NLS-equation}
\left \{
 \aligned
&\text{i}\partial_t\phi({\bf x},t)+\mathcal L \phi({\bf x},t)+\beta|\phi({\bf x},t)|^2\phi({\bf x},t)=0,\ ({\bf x},t)\in\Omega,\ t>0,\\
&\phi({\bf x},0)=\phi_0({\bf x}),\ {\bf x}\in\Omega \subset\mathbb{R}^d,\ d=1,2,3,
\endaligned
 \right.
  \end{align}
where $t$ is the time variable, ${\bf x}$ is the spatial variable, $\text{i}=\sqrt{-1}$ is the complex unit, $\phi:=\phi({\bf x},t)$ is the complex-valued wave function, $\mathcal L$ is the usual Laplace operator, $\beta$ is a given real constant, and $\phi_0({\bf x})$ is a given initial data. The nonlinear Schr\"odinger equation \eqref{NLS-equation} satisfies the following two invariants, that is mass
\begin{align}\label{NLS-mass-energy}
\mathcal M(t)=\int_{\Omega}|\phi|^2d{\bf x}\equiv \mathcal M(0),\ t\ge0,
\end{align}
and Hamiltonian energy
\begin{align}\label{NLS-Hamiltonian-energy}
\mathcal E(t) =\int_{\Omega}\Big(-|\nabla\phi|^2+\frac{\beta}{2}|\phi|^4\Big)d{\bf x}\equiv \mathcal E(0),\ t\ge0.
\end{align}

A numerical scheme that preserves one or more invariants of the original systems is called an energy-preserving method or integral-preserving method. The earlier attempts to develop energy-preserving methods for NLSE can date back to 1981, when Delfour, Fortin, and Payre \cite{DFP1981jcp} constructed a two-level finite difference scheme (also called Crank-Nicolson finite difference (CNFD) scheme) which can satisfy the discrete analogue of conservation laws \eqref{NLS-mass-energy} and \eqref{NLS-Hamiltonian-energy}. Some comparisons with nonconservation schemes are investigated by Sanz-Serna and Verwer \cite{SV1986IMA}, which show the nonconservative schemes may easily show nonlinear blow-up. Furthermore, it can be proven rigorously in mathematics that the CNFD scheme is second-order accurate in time and space \cite{BC12,BC2013KRM,Sun2005}. Up to now, the energy-preserving Crank-Nicolson schemes for NLSE are extensively extended and analyzed \cite{ADK91,ABB13,CJS1999jcp,FLM2019cicp,GWWC17,LMW2018jcam,WGX13}. Furthermore, we note that the discrete conservative laws play a crucial role in numerical analyses for numerical schemes of NLSE. However, it is fully implicit and at every time step, a nonlinear equation shall be solved by using a nonlinear iterative method and thus it may be time consuming. Other fully implicit energy-preserving schemes for NLSE provided by the averaged vector flied method \cite{QM08} and the discrete variational derivative method \cite{DO11,Furihata99} can be founded in \cite{CGM12,GCW14b,Kong2019mmas,MF01}. Zhang et al. \cite{ZPV1995amc} firstly proposed a linearly implicit Crank-Nicolson finite difference (LI-CNFD) scheme in which a linear system is to be solved at every time step. Thus it is computationally much cheaper than that of the CNFD scheme. The LI-CNFD scheme can satisfy a new discrete analogue of conservation laws \eqref{NLS-mass-energy} and \eqref{NLS-Hamiltonian-energy} and its stability and convergence is analyzed in \cite{Sun2005}. Another popular linearly implicit energy-preserving method is the relaxation finite difference scheme \cite{Besse2020IMA}. More recently, the scalar auxiliary variable (SAV) approach \cite{SXY18,SXY17} have been  a particular powerful tool for the design of linearly implicit energy-preserving numerical schemes for NLSE \cite{AST2019P}. Nevertheless, to our best knowledge, all of the schemes can achieve at most second-order in time, which cannot provide long time accurate solutions for a given large time step.

  It has been proven that no Runge-Kutta (RK) method can preserve arbitrary polynomial invariants of degree 3 or higher of arbitrary vector fields \cite{CIZ97}. Thus, over that last decades, how to develop high-order energy-preserving methods for general conservation systems attracts a lot of attention. The notable ones include, but are not limited to, high-order averaged vector field (AVF) methods \cite{LW16b,QM08}, Hamiltonian Boundary Value Methods (HBVMs) \cite{BI16,BIT10} and energy-preserving continuous stage Runge-Kutta methods  (CSRKs) \cite{CH11bit,H10,MB16,TS12}. Actually, the HBVMs and CSRKs have been shown to be an efficient method to develop high-order energy-preserving schemes for NLSE \cite{BBCIamc18,LW15}, however, the proposed schemes are fully implicit. At every time step, one needs to solve a large fully nonlinear system and thus it might be lead to high computational costs.

Due to the high computational cost of the high-order fully implicit schemes, in the literatures, ones are devoted to construct energy-preserving explicit schemes for NLSE. Based on the invariant energy quadratization (IEQ) approach \cite{YZW17} or the SAV approach and the projection method \cite{CHMR06,ELW06,Hairer00bi}, the authors proposed some explicit energy-preserving schemes for NLSE \cite{JWG2021jcam,ZhangQYS19}. Such scheme is extremely easy to implement, however, it is conditionally stable so that it may require very small time step sizes due to the certain time step restriction. Especially, this restriction is more serious in 2D and 3D. { Compared with the fully explicit scheme, the linearly implicit scheme is unconditionally stable so that the time step restriction can be removed. Furthermore, the linearly implicit scheme only need to solve a linear system of equations at every time step, thus it is more efficient than the fully implicit scheme.} In light of this, Li et al. proposed a class of linearly implicit schemes for NLSE, which can preserve the discrete analogue of conservation law \eqref{NLS-mass-energy} \cite{LGZ2020aml}. It involves solving linear equations with complicated
variable coefficients  at every time step and thus it is may be very time consuming. Another way to achieve this goal is to combine the SAV approach with the extrapolation technique or the prediction-correction strategy \cite{BWC2020arixv}. However, at every time step, the constant coefficient matrix of the scheme may admit a large norm arising from a space discretization of the linear unbounded differential operator $\mathcal L$ \cite{CCO08fcm}.

It is well-known that the exponential integrator often involve exact integration of the linear part of the target equation, and thus it can achieve high accuracy, stability for a very stiff differential equation, such as highly oscillatory ODEs and semidiscrete time-dependent PDEs. For more details on the exponential integrator, please refer to the excellent review article provided by Hochbruck and Ostermann \cite{HO10an}. As a matter of fact, many energy-preserving exponential integrators have been done for the conservative systems in the past few decades. Based on the projection approach, Celledoni et al. developed some  symmetry- and energy-preserving implicit exponential schemes for the cubic Schr\"odinger equation \cite{CCO08fcm}. By combining the exponential integrator with the AVF method, Li and Wu constructed a class of second-order implicit energy-preserving exponential schemes for  canonical Hamiltonian systems and successfully applied it to solve NLES \cite{LW16b}. Further analysis and generalization is investigated by Shen and Leok \cite{SL19jcp}. More recently, Jiang et al. showed that the SAV approach is also an efficient approach to develop second-order linearly implicit energy-preserving exponential integrator for NLSE \cite{JWC2020jcp}. Overall, there exist very few works devoted to development of high-order linearly implicit energy-preserving exponential schemes for NLSE, which motivates this paper. We should note that, in general, the particularly interesting types
of exponential integrators are integrating factor (IF) methods and exponential time differencing (ETD) methods,
respectively and we here mainly focus on the high-order linearly implicit energy-preserving IF schemes.

In this paper, following the idea of the SAV approach, we firstly reformulate the original system \eqref{NLS-equation} into a new system by introducing a new auxiliary variable, which satisfies a modified energy. The spatial discretization is then performed with the standard Fourier pseudo-spectral method \cite{CQ01,ST06}. Subsequently, the extrapolation technique/ prediction-correction strategy is employed to the nonlinear terms of the semi-discrete system and a linearized energy-conserving system is obtained. Based on the Lawson transformation, the semi-discrete system is rewritten as an equivalent form and the fully discrete scheme is obtained by using a RK method to the linearized equation. We show that, the proposed scheme can preserve the modified energy in the discretized level when some under certain circumstances are imposed on the coefficients of the Runge-Kutta method and at every time step, only a linear equations with constant coefficients is solved in which there is no existing references \cite{BWC2020arixv,LS2020jsc} considering this issue.

The rest of this paper is organized as follows. In Section \ref{Sec:LI-EP:2}, we use the idea of the SAV approach to reformulate the original equation \eqref{NLS-equation} into a new reformulation, which satisfies a modified energy. In Section \ref{Sec:LI-EP:3}, two class of linearly implicit energy-preserving integrating
factor Runge-Kutta methods are proposed for NLSE \eqref{NLS-equation}. In Section \ref{Sec:LI-EP:4}, several numerical examples are investigated to illustrate the efficiency of the proposed scheme. We draw some conclusions in Section \ref{Sec:LI-EP:5}.

\section{Model reformulation}\label{Sec:LI-EP:2}
In this section, we first employ the SAV approach to reformulate the NLSE \eqref{NLS-equation} into a new reformulation which satisfies a quadratic energy. The resulting reformulation provides an elegant platform for developing high-order linearly implicit energy-preserving schemes.

Based on the idea of the SAV approach \cite{SXY18,SXY17}, we introduce an auxiliary variable
\begin{align*}
p:=p(t)=\sqrt{(\phi^2,\phi^2)+C_0},
\end{align*}
where  $C_0\ge 0$ is a constant to make $p$ well-defined for all $\phi$. Here, $(v,w)$ is the continuous inner product defined by $\int_{\Omega}v\bar{w}d{\bf x}$ where $\bar{w}$ represents the conjugate of $w$. Then, the Hamiltonian energy functional is rewritten as the following quadratic form
\begin{align*}
\mathcal E(t)=( {\mathcal L}\phi, \phi)+\frac{\beta}{2}p^2-\frac{\beta}{2}C_0.
\end{align*}
According to the energy variational principle, we obtain the following SAV reformulated system
\begin{align}\label{NLS-SAV-reformulation}
\left\{
\begin{aligned}
&\partial_t \phi=\text{i}\bigg ( {\mathcal L}\phi+ \frac{\beta|\phi|^2\phi p}{\sqrt{(\phi^2,\phi^2)+C_0}} \bigg), \\
&\frac{d}{dt} p=2 {\rm Re}\bigg(\text{i} {\mathcal L}\phi,\frac{|\phi|^2\phi}{\sqrt{(\phi^2,\phi^2)+C_0}} \bigg),
\end{aligned}
\right.
\end{align}
with the consistent initial condition
\begin{align}\label{eq2.6}
\phi({\bf x},0)=\phi_0({\bf x}),  \  p(0)=\sqrt{\big(  \phi^2_0({\bf x}), \phi^2_0({\bf x})) +C_0   },
\end{align}
where ${\rm Re}(u)$ represents the real part of $u$.
 \begin{thm}\label{LI-E-NLS-thm-2.1} The SAV reformulation \eqref{NLS-SAV-reformulation} preserves the following quadratic energy
  \begin{align}\label{Q-energy-conservetion-law}
  \frac{d}{dt} \mathcal{E}(t)=0,\ \mathcal{E}(t)=(\mathcal{L}\phi,\phi)+\frac{\beta}{2}p^2-\frac{\beta}{2}C_0,\ t\ge0.
  \end{align}
  \end{thm}
  \begin{prf} It is clear to see
  \begin{align*}
  \frac{d}{dt} \mathcal{E}(t)&=2\text{Re}(\mathcal{L}\phi,\partial_t\phi)+\beta p\partial_tp\\
  &=2\beta p\text{Re}\Big(\mathcal{L}\phi,\text{\rm i}\frac{|\phi|^2\phi}{\sqrt{(\phi^2,\phi^2)+C_0}}\Big)+\beta p\partial_tp\\
  &=-\beta p\partial_t p+\beta p\partial_tp\\
  &=0.
  \end{align*}
  This completes the proof. \qed
  \end{prf}

\section{The construction of the high-order linearly integrating
factor Runge-Kutta methods}\label{Sec:LI-EP:3}
In this section, a class of high-order linearly implicit energy-preserving integrating
factor Runge-Kutta methods are proposed for the SAV reformulated system \eqref{NLS-SAV-reformulation}. For simplicity, in this paper, we shall introduce our schemes in two space dimension, i.e., $d=2$ in \eqref{NLS-SAV-reformulation}. Generalizations to $d=1$ or $3$ are straightforward.

 \subsection{Spatial discretization}\label{Sec:PM:3}
  We set the computational domain $\Omega=[a,b]\times[c,d]$ and let $l_x=b-a$ and $l_y=d-c$. Choose the mesh sizes $h_x=l_x/N_x$ and $h_y=l_y/N_y$ with two even positive integers $N_x$ and $N_y$, and denote the grid points by $x_{j}=jh_x$ for $j=0,1,2,\cdots,N_x$ and $y_{k}=kh_y$ for $k=0,1,2,\cdots,N_y$; let $\phi_{j,k}$ be the numerical approximation of  $\phi(x_j,y_k,t)$ for $j=0,1,\cdots,N_x,k=0,1,2,\cdots,N_y$, and
  $$\phi:=(\phi_{0,0},\phi_{1,0},\cdots,\phi_{N_x-1,0},\phi_{0,1},\phi_{1,1},\cdots,\phi_{N_x-1,1},\cdots,\phi_{0,N_y-1},\phi_{1,N_y-1},\cdots,\phi_{N_x-1,N_y-1})^{T}$$ be the solution vector; we also define the discrete inner product, $L^2$-norm and $L^{\infty}$-norm as, respectively,
\begin{align*}
&\langle \phi,\psi\rangle_{h}=h_xh_y\sum_{j=0}^{N_x-1}\sum_{k=0}^{N_y-1}\phi_{j,k}\bar{\psi}_{j,k},\ \|\phi\|_{h}^2=\langle \phi,\phi\rangle_{h},\ \|{ \phi}\|_{h,\infty}=\max\limits_{j,k}|\phi_{j,k}|.
\end{align*}
  In addition, we denote $`\cdot$' as the element product of vectors ${\phi}$ and  ${\psi}$, that is
\begin{align*}
\phi\cdot \psi=&\big(\phi_{0,0}\psi_{0,0},\cdots,\phi_{N_x-1,0}\psi_{N_x-1,0},\cdots,\phi_{0,N_x-1}\psi_{0,N_y-1},\cdots,\phi_{N_x-1,N_y-1}\psi_{N_x-1,N_y-1}\big)^{T}.
\end{align*}
For brevity, we denote ${\phi}\cdot\phi$ and ${\phi}\cdot\bar{\phi}$  as ${\phi}^2$ and $|\phi|^2$, respectively.

Denote the interpolation space as
\begin{align*}
&\mathcal T_h=\text{span}\{g_{j}(x)g_{k}(y),\ 0\leq j\leq N_x-1, 0\leq k\leq N_y-1\}
\end{align*}
where $g_{j}(x)$ and $g_{k}(y)$ are trigonometric polynomials of degree $N_{x}/2$ and $N_{y}/2$,
given, respectively, by
\begin{align*}
  &g_{j}(x)=\frac{1}{N_{x}}\sum_{l=-N_{x}/2}^{N_{x}/2}\frac{1}{a_{l}}e^{\text{i}l\mu_{x} (x-x_{j})},\ g_{k}(y)=\frac{1}{N_{y}}\sum_{l=-N_{y}/2}^{N_{y}/2}\frac{1}{b_{l}}e^{\text{i}l\mu_{y} (y-y_{k})},
\end{align*}
with $a_{l}=\left \{
 \aligned
 &1,\ |l|<\frac{N_x}{2},\\
 &2,\ |l|=\frac{N_x}{2},
 \endaligned
 \right.,\ b_{l}=\left \{
 \aligned
 &1,\ |l|<\frac{N_y}{2},\\
 &2,\ |l|=\frac{N_y}{2},
 \endaligned
 \right.$, $\mu_{x}=\frac{2\pi}{l_{x}}$ and $\mu_{y}=\frac{2\pi}{l_{y}}$.
We then define the interpolation operator $I_{N}: C(\Omega)\to \mathcal T_h$ as:
\begin{align*}
I_{N}\phi(x,y,t)=\sum_{j=0}^{N_{x}-1}\sum_{k=0}^{N_{y}-1}\phi_{j,k}(t)g_{j}(x)g_{k}(y),
\end{align*}
where $\phi_{j,k}(t)=\phi(x_{j},y_{k},t)$. 

Taking the derivative with respect to $x$ and $y$, respectively, and then evaluating the resulting expressions at the collocation points ($x_{j},y_{k}$), we have
\begin{align}\label{GR-KDV-2.2}
\frac{\partial^{2} I_{N}\phi(x_{j},y_{k},t)}{\partial x^{2}}
&=\sum_{l=0}^{N_{x}-1}\phi_{l,k}(t)\frac{d^{2}g_{l}(x_{j})}{dx^{2}}=\sum_{l=0}^{N_{x}-1}(D_2^{x})_{j,l}\phi_{l,k}(t),\\
\frac{\partial^{2} I_{N}\phi(x_{j},y_{k},t)}{\partial y^{2}}&=\sum_{l=0}^{N_{y}-1}\phi_{j,l}\frac{d^{2}g_{l}(y_{k})}{dy^{2}}=\sum_{l=0}^{N_{y}-1}\phi_{j,l}(t)(D_2^y)_{k,l},
\end{align}
where
\begin{align*}
 &(D^{x}_{2})_{j,k}=
 \left \{
 \aligned
 &\frac{1}{2}\mu_{x}^{2} (-1)^{j+k+1}\csc^{2}(\mu_{x} \frac{x_{j}-x_{k}}{2}),\ &j\neq k,\\
 &-\mu_{x}^{2}\frac{N_{x}^{2}+2}{12},\quad \quad \quad \quad ~ &j=k,
 \endaligned
 \right.\\
 &(D^{y}_{2})_{j,k}=
 \left \{
 \aligned
 &\frac{1}{2}\mu_{y}^{2} (-1)^{j+k+1}\csc^{2}(\mu_{y} \frac{y_{j}-y_{k}}{2}),\ &j\neq k,\\
 &-\mu_{y}^{2}\frac{N_{y}^{2}+2}{12},\quad \quad \quad \quad ~ &j=k.
 \endaligned
 \right.
 \end{align*}
   \begin{rmk} We should note that \cite{ST06}
   \begin{align*}
D_2^w=\mathcal F_w^H\Lambda_w \mathcal F_w,\ w=x,y,
\end{align*}
where
\begin{align*}
 \Lambda_w
=\Big(\text{\rm i}\frac{2\pi}{l_w}\Big)^2\text{\rm diag}\Big[0,1,\cdots,\frac{N_w}{2},-\frac{N_w}{2}+1,\cdots,-2,-1\Big]^2,
\end{align*}
and $\mathcal F_w$ is the discrete Fourier transform (DFT) and $\mathcal F_w^H$ represents the conjugate
transpose of $\mathcal F_w$.
   \end{rmk}

  Then, we use the standard Fourier pseudo-spectral method to solve \eqref{NLS-SAV-reformulation} and obtain
    \begin{align}\label{NLS-SP-SAV-reformulation}
\left\{
\begin{aligned}
&\frac{d}{dt} \phi=\text{i}\bigg ( {\mathcal L}_h\phi+ \frac{\beta|\phi|^2\cdot \phi p}{\sqrt{\langle\phi^2,\phi^2\rangle_h+C_0}} \bigg), \\
&\frac{d}{dt} p=2 {\rm Re}\bigg\langle\text{i} {\mathcal L}_h\phi,\frac{|\phi|^2\cdot\phi}{\sqrt{\langle\phi^2,\phi^2\rangle_h+C_0}} \bigg\rangle_h,
\end{aligned}
\right.
\end{align}
where ${\mathcal L}_h=I_{N_y}\otimes D_2^x+D_2^y\otimes I_{N_x}$ is the spectral differentiation matrix and $\otimes$ represents the Kronecker product.
   \begin{thm} The semi-discrete system \eqref{NLS-SP-SAV-reformulation} preserves the following semi-discrete quadratic energy
  \begin{align}\label{sm-SAV-energy-conservetion-law}
  \frac{d}{dt} E_h(t)=0,\ E_h(t)=\langle\mathcal{L}_h\phi,\phi\rangle_h+\frac{\beta}{2}p^2-\frac{\beta}{2}C_0,\ t\ge0.
  \end{align}
   \end{thm}
   \begin{prf} The proof is similar to Theorem \ref{LI-E-NLS-thm-2.1}, thus, for brevity, we omit it.

   \end{prf}

   \subsection{Temporal exponential integration}\label{Sec:PM:3}
   Denote $t_{n}=n\tau,$ and $t_{ni}=t_n+c_i\tau,\ i=1,2,\cdots,s$, $n = 0,1,2,\cdots,$ where $\tau$ is the time step and $c_1,c_2,\cdots,c_s$ are distinct real numerbers (usually $0\le c_i\le 1$). The approximations of the function $\phi(x,y,t)$ at points  $(x_j,y_k,t_n)$ and $(x_j,y_k,t_{ni})$ are denoted by $\phi_{j,k}^n$ and $(\phi_{ni})_{j,k}$, and the approximations of the function $p(t)$ at points  $t_n$ and $t_{ni}$ are denoted by $p^n$ and $p_{ni}$.
\subsubsection{ Linearly implicit integrating
factor Runge-Kutta method based on the extrapolation method}
{We assume that the numerical solutions of $\phi$ for \eqref{NLS-SP-SAV-reformulation} from $t=0$ to $t\le t_n$ are obtained. Then we obtain the Lagrange interpolating
polynomial (denoted by $\mathcal P(t)$) based on the time nodes $t_k,\ t_{k}+c_i\tau,$ and the values $\phi^k,\ \phi_{ki}$, where $\ 0\le k\le n,\ i=1,2,\cdots,s$. Subsequently, we use $\mathcal P(t)$ to approximate $\phi$ of the nonlinear terms of \eqref{NLS-SP-SAV-reformulation} over time interval   $(t_n,t_{n+1}]$ and a linearized system is obtained, as follows:
\begin{align}\label{linear-NLS-SP-SAV-reformulation}
\left\{
\begin{aligned}
&\frac{d}{dt} \tilde\phi=\text{i}\bigg ( {\mathcal L}_h\tilde\phi+ \frac{\beta|\mathcal P(t)|^2\cdot \mathcal P(t) \tilde p}{\sqrt{\langle\mathcal P(t)^2,\mathcal P(t)^2\rangle_h+C_0}} \bigg), \\
&\frac{d}{dt} \tilde p=2 {\rm Re}\bigg\langle\text{i} {\mathcal L}_h\tilde\phi,\frac{|\mathcal P(t)|^2\cdot\mathcal P(t)}{\sqrt{\langle\mathcal P(t)^2,\mathcal P(t)^2\rangle_h+C_0}} \bigg\rangle_h,\ t\in(t_n,t_{n+1}],
\end{aligned}
\right.
\end{align}}
where $n=0,1,2,\cdots,$ $\tilde\phi:=\tilde\phi(t)$ and $\tilde p:=\tilde p(t)$ are approximations of $\phi(t)$ and $p(t)$, respectively over time interval $(t_n,t_{n+1}],\ n=0,1,2,\cdots,$.

\begin{rmk} It is remarked that, by the similar argument as Theorem \ref{LI-E-NLS-thm-2.1}, we can show the linearized system \eqref{linear-NLS-SP-SAV-reformulation} satisfies a semi-discrete
quadratic energy, as follows:
 \begin{align*}
  \frac{d}{dt} \tilde E_h(t)=0,\ \tilde E_h(t)=\langle\mathcal{L}_h\tilde\phi,\tilde\phi\rangle_h+\frac{\beta}{2}(\tilde p)^2-\frac{\beta}{2}C_0,\ t\ge0.
  \end{align*}
\end{rmk}


By using the Lawson transformation~\cite{Lawson1967}, we multiply both sides of the first equation of~\eqref{linear-NLS-SP-SAV-reformulation} by the operator $\rm exp(-{\rm i}{\mathcal L}_h t)$, and then introduce $\Psi=\rm exp(-{\rm i}{\mathcal L}_h t){\color{blue}\tilde\phi}$ to transform~\eqref{linear-NLS-SP-SAV-reformulation} into an equivalent form, as follows:
\begin{align}\label{Lawson-NLS-SP-SAV-reformulation}
\left\{
\begin{aligned}
&\frac{d}{dt}\Psi={\rm i}\exp(-{\rm i}\mathcal{L}_ht)\frac{\beta|\mathcal P(t)|^2\cdot\mathcal P(t) {\color{blue}\tilde p}}{\sqrt{\langle\mathcal P(t)^2,\mathcal P(t)^2\rangle_h+C_0}}, \\
&\frac{d}{dt} {\color{blue} \tilde p} =\!2{\rm Re}\Bigg\langle\!\!{\rm i} \exp({\rm i}\mathcal{L}_ht) {\mathcal L}_h\Psi,\frac{|\mathcal P(t)|^2\cdot\mathcal P(t)}{\sqrt{\big\langle\mathcal P(t)^2,\mathcal P(t)^2\big\rangle_h+C_0}}  \Bigg\rangle_h,\ t\in(t_n,t_{n+1}],
\end{aligned}
\right.
\end{align}
where $n=0,1,2,\cdots,$.

 We first apply an RK method to the linearized system \eqref{Lawson-NLS-SP-SAV-reformulation}, and then the discretization is rewritten in terms of the original variables to give a class of linearly implicit integrating
factor Runge-Kutta methods (LI-IF-RKMs), as follows:
\begin{shm}\label{NLS-shn-3.1}
 Let $b_i,a_{ij}\ (i,j=1,\cdots,s)$ be real numbers and let $c_i=\sum_{j=1}^sa_{ij}$. For given $\phi^n$, $p^n$ and $\phi_{ni}^{*}=\mathcal P(t_n+c_i\tau),i=1,\cdots,s$, $s$-stage linearly implicit integrating factor Runge-Kutta method is given by
\begin{align}\label{LI-E-scheme1}
\left\{
\begin{aligned}
&\phi_{ni}=\exp({\rm i}\mathcal{L}_hc_i\tau) \phi^n+\tau\sum_{j=1}^{s}a_{ij}\exp({\rm i}\mathcal{L}_h(c_i-c_j)\tau) k_j,\\
&p_{ni}=p^n+\tau\sum_{j=1}^{s}a_{ij}l_j,\ l_i\!=2{\rm Re}\!\Big\langle\!{\rm i}\mathcal{L}_h\phi_{ni}, \frac{|\phi_{ni}^{*}|^2\cdot\phi_{ni}^{*}}
{\sqrt{\big\langle(\phi_{ni}^{*})^2,(\phi_{ni}^{*})^2\big\rangle_h\!+\!C_0}}\Big\rangle_h,\\
&k_i\!=\!{\rm i}\ \frac{\beta|\phi_{ni}^{*}|^2\cdot\phi_{ni}^{*} p_{ni}}
{\sqrt{\big\langle(\phi_{ni}^{*})^2,(\phi_{ni}^{*})^2\big\rangle_h\!+\!C_0}},
\end{aligned}
\right.
\end{align}
where $i=1,2,\cdots,s$.
Then $(\phi^{n+1},p^{n+1})$ are updated by
\begin{align}\label{LI-E-scheme2}
\left\{
\begin{aligned}
&\phi^{n+1}=\exp({\rm i}\mathcal{L}_h\tau) \phi^n+\tau\sum_{i=1}^{s}b_i \exp({\rm i}\mathcal{L}_h(1-c_i)\tau) k_i,\\
&p^{n+1}=p^n+\tau\sum_{i=1}^sb_il_i.
\end{aligned}
\right.
\end{align}
\end{shm}

\begin{PPS}\label{NLS-exp-matrix}  According to the definition of the matrix-valued function \cite{Higham2008siam}, it holds
\begin{itemize}
\item $\exp({\rm i}\mathcal{L}_h\tau)^{H}=\exp(-{\rm i}\mathcal{L}_h\tau)$ and $\exp({\rm i}\mathcal{L}_h\tau)\mathcal L_h=\mathcal L_h\exp({\rm i}\mathcal{L}_h\tau)$;
 \item $\exp({\rm i}\mathcal{L}_h\tau)=\mathcal F_{N_y}^{H}\otimes\mathcal F_{N_x}^{H}\exp({\rm i}(I_{N_y}\otimes\Lambda_x+\Lambda_y\otimes I_{N_x})\tau)\mathcal F_{N_y}\otimes\mathcal F_{N_x}$, which can be efficiently implemented by the matlab functions \emph{fftn.m} and \emph{ifftn.m}.
\end{itemize}
\end{PPS}

\begin{thm} \label{th2.1}  If the RK coefficients of { \eqref{LI-E-scheme1} and \eqref{LI-E-scheme2}} satisfy
\begin{align}\label{LI-EI-coefficient}
b_ia_{ij}+b_ja_{ji}=b_ib_j, \ \forall \ i,j=1,\cdots,s,
\end{align}
 the proposed scheme \eqref{LI-E-scheme1}-\eqref{LI-E-scheme2} can preserve the discrete quadratic energy, as follows:
\begin{align}\label{LI-EI-F-energy}
E_h^{n}=E_h^0, \ E_h^n=\langle\mathcal{L}_h\phi^n,\phi^n\rangle_h+\frac{\beta}{2}(p^n)^2-\frac{\beta}{2}C_0,\  n=0,1,2,\cdots.
\end{align}
\end{thm}
\begin{prf}
We obtain from the first equality of~\eqref{LI-E-scheme2}, Proposition \ref{NLS-exp-matrix} and \eqref{LI-EI-coefficient} that
\begin{align*}
&\langle{\mathcal L}_h\phi^{n+1},\phi^{n+1} \rangle_h-\langle{\mathcal L}_h\phi^{n},\phi^{n}\rangle_h \nonumber\\
=&2\tau\sum_{i=1}^{s}b_i {\rm Re}\big\langle \exp({\rm i}\mathcal{L}_hc_i\tau)\phi^n,\mathcal{L}_hk_i \big\rangle_h+\tau^2 \sum_{i,j=1}^{s} b_ib_j \big\langle\mathcal{L}_hk_i,\exp({\rm i}\mathcal{L}_h(c_i-c_j)\tau) k_j \big \rangle_h\nonumber\\
=& 2\tau\sum_{i=1}^{s}b_i {\rm Re}\big\langle  \exp({\rm i}\mathcal{L}_hc_i\tau)\phi^n,\mathcal L_hk_i \big\rangle_h
  +\tau^2 \sum_{i,j=1}^{s} (b_ia_{ij}+b_ja_{ji}) \big\langle \mathcal{L}_h k_i,\exp({\rm i}\mathcal{L}_h(c_i-c_j)\tau) k_j \big \rangle_h\\
=& 2\tau\sum_{i=1}^{s}b_i {\rm Re}\big\langle  \exp({\rm i}\mathcal{L}_hc_i\tau)\phi^n,\mathcal L_hk_i \big\rangle_h
  +2\tau^2 \sum_{i,j=1}^{s} b_ia_{ij}{\rm Re}\big\langle\exp({\rm i}\mathcal{L}_h(c_i-c_j)\tau) k_j, \mathcal{L}_h k_i\big \rangle_h.
  \end{align*}
  With the first equality of~\eqref{LI-E-scheme1}, we have
  \begin{align*}
  2\tau\sum_{i=1}^{s}b_i{\rm Re}&\big\langle  \exp({\rm i}\mathcal{L}_hc_i\tau)\phi^n,\mathcal L_hk_i \big\rangle_h\nonumber\\
  &=2\tau\sum_{i=1}^{s}b_i{\rm Re}\big\langle  \phi_{ni},\mathcal L_hk_i \big\rangle_h-2\tau^2\sum_{i,j=1}^{s}b_ia_{ij}{\rm Re}\big\langle\exp({\rm i}\mathcal{L}_h(c_i-c_j)\tau) k_j,\mathcal L_hk_i \big\rangle_h.
  \end{align*}
  Thus, using above two equalities, we can obtain
\begin{align}\label{LI-EI-NLS-3.16}
&\langle{\mathcal L}_h\phi^{n+1},\phi^{n+1} \rangle_h-\langle{\mathcal L}_h\phi^{n},\phi^{n}\rangle_h=2\tau \sum_{i=1}^s b_i{\rm Re}\big\langle\!\phi_{ni},\mathcal L_hk_i \big\rangle_h.
\end{align}
Moreover, we derive from \eqref{LI-E-scheme1}, the second equality of \eqref{LI-E-scheme2} and \eqref{LI-EI-coefficient} that
\begin{align}\label{LI-EI-NLS-3.17}
\frac{\beta}{2}\big[(p^{n+1})^2-(p^n)^2\big]
&=\beta\tau \sum_{i=1}^s b_il_i p^n+\frac{\beta}{2}\tau^2\sum_{i,j=1}^sb_ib_jl_il_j\nonumber\\
&= \beta\tau \sum_{i=1}^s b_il_i \big(p_{ni}-\tau\sum_{j=1}^sa_{ij}l_j\big)+\frac{\beta}{2}\tau^2\sum_{i,j=1}^s(b_ia_{ij}+b_ja_{ji})l_il_j  \nonumber\\
&=2\tau \sum_{i=1}^s b_ip_{ni}{\rm Re}\bigg\langle\!{\rm i}\mathcal L_h\phi_{ni},\frac{\beta|\phi_{ni}^{*}|^2\cdot\phi_{ni}^{*}}{\sqrt{\langle(\phi_{ni}^{*})^2,(\phi_{ni}^{*})^2\rangle_h+C_0}}\bigg\rangle_h\nonumber\\
&= -2\tau \sum_{i=1}^s b_i{\rm Re}\big\langle\!\phi_{ni},\mathcal L_hk_i \big\rangle_h.
\end{align}
It follows from~\eqref{LI-EI-NLS-3.16} and~\eqref{LI-EI-NLS-3.17}, we completes the proof.\qed
\end{prf}
\begin{rmk} { The Gauss collocation method where the RK coefficients $c_1,c_2,\cdots,c_s$ are chosen as the Gaussian quadrature
nodes, i.e., the zeros of the $s$-th shifted Legendre polynomial $\frac{d^s}{dx^s}(x^s(x-1)^s)$ satisfies the condition \eqref{LI-EI-coefficient} (see \cite{ELW06}). In addition, it is proven that the $s$-stage Gauss collocation method has the same order $2s$ as the Gaussian quadrature formula. In particular, the RK coefficients for 2-stage and 3-stage Gauss collocation methods can be given explicitly by (see \cite{ELW06,Sanzs88})}

\begin{table}[H]
\centering
\begin{tabular}{c|cc}
$\frac{1}{2}-\frac{\sqrt{3}}{6}$ &$\frac{1}{4}$ & $\frac{1}{4}- \frac{\sqrt{3}}{6}$\\
$\frac{1}{2}+\frac{\sqrt{3}}{6}$ &$\frac{1}{4}+ \frac{\sqrt{3}}{6}$ &$\frac{1}{4}$ \\
\hline
                                 &$\frac{1}{2}$&$\frac{1}{2}$
\end{tabular},\ \ \
\begin{tabular}{c|ccc}
$\frac{1}{2}-\frac{\sqrt{15}}{10}$ &$\frac{5}{36}$ &  $\frac{2}{9}-\frac{\sqrt{15}}{15}$    &$\frac{5}{36}- \frac{\sqrt{15}}{30}$\\
$\frac{1}{2}$   &$\frac{5}{36}+ \frac{\sqrt{15}}{24}$  &  $\frac{2}{9}$  & $\frac{5}{36}- \frac{\sqrt{15}}{24}$ \\
$\frac{1}{2}+\frac{\sqrt{15}}{10}$ & $\frac{5}{36}+ \frac{\sqrt{15}}{30}$ & $\frac{2}{9}+\frac{\sqrt{15}}{15}$  & $\frac{5}{36}$ \\
\hline
& $\frac{5}{18}$ & $\frac{4}{9}$  & $\frac{5}{18}$
\end{tabular}.
\caption{The RK coefficients of 2-stage (left) and 3-stage (right) Gauss collocation methods.}\label{LI-EI-tab-GM}
\end{table}
\end{rmk}

\begin{rmk}\label{NLS-rmk-3.4}  It is well-known that the interpolation polynomial will be highly oscillating when too 
many interpolation points are chosen, which makes $\phi_{ni}^*$ be an inaccurate extrapolation for $\phi(t_{ni})$. Thus, in this paper, {we only choose $t_{n-1}$, $t_{n-1}+c_i\tau,\ i=1,2,\cdots,s$ and $t_{n}$ as the interpolation points, together with the values $\phi^{n-1},\ \phi_{(n-1)i},\ i=1,2,\cdots,s$ and $\phi^n$, for the interpolation, where the coefficients $c_1,c_2,\cdots,c_s$ are chosen as the Gaussian quadrature
nodes. The resulting extrapolation $\phi_{ni}^*$ may provide the approximation of the order $\mathcal O(\tau^{s+1})$ for $\phi(t_{ni})$ \cite{LS2020jsc}. Thus, if such extrapolation and the RK coefficients of the $s$-stage Gauss collocation method are chosen for \eqref{LI-E-scheme1}, the obtained scheme may achieve ($s+1$)-order accuracy in time.} For example, we  choose $(t_{n-1},\phi^{n-1}),\ (t_{n-1}+c_1\tau,\phi_{(n-1)1}),\ (t_{n-1}+c_2\tau,\phi_{(n-1)2})$ and $(t_{n},\phi^{n})$, where $c_1=\frac{1}{2}-\frac{\sqrt{3}}{6}$ and $c_2=\frac{1}{2}+\frac{\sqrt{3}}{6}$, the Lagrange interpolating polynomial at $t_n+c_1\tau$ and $t_n+c_2\tau$ is \cite{GZW2020jcp}
\begin{align}\label{2-stage-EI-NLS1}
&\phi_{n1}^{*} = (2\sqrt3-4)\phi^{n-1}+(7\sqrt3-11)\phi_{(n-1)1}\nonumber\\
&~~~~~~~~~~~~~~~~~~~~~~~~~~~~~~~~~~~~~~~~+(6-5\sqrt3)\phi_{(n-1)2}+(10-4\sqrt3)\phi^{n},\\\label{2-stage-EI-NLS2}
&\phi_{n2}^{*} = -(2\sqrt3+4)\phi^{n-1}+(6+5\sqrt3)\phi_{(n-1)1}\nonumber\\
&~~~~~~~~~~~~~~~~~~~~~~~~~~~~~~~~~~~~~~~-(7\sqrt3+11)\phi_{(n-1)2}+(10+4\sqrt3)\phi^n.
\end{align}
Then, the scheme \eqref{LI-E-scheme1} based on the extrapolation \eqref{2-stage-EI-NLS1}-\eqref{2-stage-EI-NLS2} and the RK coefficients of the 2-stage Gauss method (see Table \ref{LI-EI-tab-GM}) may achieve third-order accuracy in time. In addition, if we choose $(t_{n-1},\phi^{n-1}),\ (t_{n-1}+c_1\tau,\phi_{(n-1)1}),\ (t_{n-1}+c_2\tau,\phi_{(n-1)2})$ and $(t_{n-1}+c_3\tau,\phi_{(n-1)3})$, where $c_1=\frac{1}{2}-\frac{\sqrt{15}}{10}$, $c_2=\frac{1}{2}$  and $c_3=\frac{1}{2}+\frac{\sqrt{15}}{10}$, the Lagrange interpolating polynomial at $t_n+c_1\tau$, $t_n+c_2\tau$ and $t_n+c_3\tau$ is given by \cite{BWC2020arixv}
\begin{align}\label{3-stage-EI-NLS1}
&\phi_{n1}^*=(6\sqrt{15}-26)\phi^{n-1}+(-5\sqrt{15}/3+11)\phi_{(n-1)1}\nonumber\\
&~~~~~~~~~~~~~~~~~~~~~~~+(16\sqrt{15}/3-24)\phi_{(n-1)2}+(-29\sqrt{15}/3+40)\phi_{(n-1)3},\\\label{3-stage-EI-NLS2}
&\phi_{n2}^*=-17\phi^{n-1}+(5\sqrt{15}/2+\frac{35}{2})\phi_{(n-1)1}-17\phi_{(n-1)2}\nonumber\\
&~~~~~~~~~~~~~~~~~~~~~~~~~~~~~~~~~~~~~~~~~~~~~~~~~~+(-5\sqrt{15}/2+\frac{35}{2})\phi_{(n-1)3},\\\label{3-stage-EI-NLS3}
&\phi_{n3}^*=(-6\sqrt{15}-26)\phi^{n-1}+(29\sqrt{15}/3+40)\phi_{(n-1)1}\nonumber\\
&~~~~~~~~~~~~~~~~~~~~~~~+(-16\sqrt{15}/3-24)\phi_{(n-1)2}+(5\sqrt{15}/3+11)\phi_{(n-1)3}.
\end{align}
 The scheme \eqref{LI-E-scheme1} based on the extrapolation \eqref{3-stage-EI-NLS1}-\eqref{3-stage-EI-NLS3} and the RK coefficients of the 3-stage Gauss method (see Table \ref{LI-EI-tab-GM}) may achieve fourth-order accuracy in time.
\end{rmk}

Besides its energy-preserving property, a most remarkable thing about the above
scheme \eqref{LI-E-scheme1}-\eqref{LI-E-scheme2} is that it can be solved efficiently. For simplicity, we take $s=2$ as an example.

We denote
\begin{align}\label{NLS-fast-solver1}
\gamma_i^{*} = \text{i}\frac{|\phi_{ni}^{*}|^2\cdot\phi_{ni}^{*}}{\sqrt{\langle(\phi_{ni}^{*})^2,(\phi_{ni}^{*})^2\rangle_h+C_0}},\ i=1,2,
\end{align}
and rewrite
\begin{align}\label{NLS-fast-solver2}
&l_1=-2\text{Re}\Big\langle\mathcal L_h\phi_{n1},\gamma_1^{*}\Big\rangle_h,\ l_2=-2\text{Re}\Big\langle\mathcal L_h\phi_{n2},\gamma_2^{*}\Big\rangle_h.
\end{align}
With \eqref{NLS-fast-solver1}, we have
\begin{align}\label{NLS-fast-solver3}
&k_1 = \beta\gamma_1^{*}p^n-2\beta\tau a_{11}\gamma_1^{*}\text{Re}\Big\langle\mathcal L_h\phi_{n1},\gamma_1^{*}\Big\rangle_h-2\beta\tau a_{12}\gamma_1^{*}\text{Re}\Big\langle\mathcal L_h\phi_{n2},\gamma_2^{*}\Big\rangle_h,\\\label{NLS-fast-solver4}
&k_2 = \beta\gamma_2^{*}p^n-2\beta\tau a_{21}\gamma_2^{*}\text{Re}\Big\langle\mathcal L_h\phi_{n1},\gamma_1^{*}\Big\rangle_h-2\beta\tau a_{22}\gamma_2^{*}\text{Re}\Big\langle\mathcal L_h\phi_{n2},\gamma_2^{*}\Big\rangle_h.
\end{align}
Then it follows from the first equality of \eqref{LI-E-scheme1} and \eqref{NLS-fast-solver3}-\eqref{NLS-fast-solver4} that
\begin{align}\label{NLS-fast-solver5}
\phi_{n1}&=\gamma_1^n-(2\beta\tau^2a_{11}^2\gamma_1^{*}+2\beta\tau^2a_{12}a_{21}\exp(\text{i}\mathcal L_h(c_1-c_2)\tau)\gamma_2^{*})\text{Re}\Big\langle\mathcal L_h\phi_{n1},\gamma_1^{*}\Big\rangle_h\nonumber\\
&-(2\beta\tau^2a_{11}a_{12}\gamma_1^{*}+2\beta\tau^2a_{12}a_{22}\exp(\text{i}\mathcal L_h(c_1-c_2)\tau)\gamma_2^{*})\text{Re}\Big\langle\mathcal L_h\phi_{n2},\gamma_2^{*}\Big\rangle_h,\\\label{NLS-fast-solver6}
\phi_{n2}&=\gamma_2^n-(2\beta\tau^2a_{21}a_{11}\exp(\text{i}\mathcal L_h(c_2-c_1)\tau)\gamma_1^{*}+2\beta\tau^2a_{22}a_{21}\gamma_2^{*})\text{Re}\Big\langle\mathcal L_h\phi_{n1},\gamma_1^{*}\Big\rangle_h\nonumber\\
&-(2\beta\tau^2a_{21}a_{12}\exp(\text{i}\mathcal L_h(c_2-c_1)\tau)\gamma_1^{*}+2\beta\tau^2a_{22}^2\gamma_2^{*})\text{Re}\Big\langle\mathcal L_h\phi_{n2},\gamma_2^{*}\Big\rangle_h,
\end{align}
where
\begin{align*}
&\gamma_1^n=\exp(\text{i}\mathcal L_h c_1\tau)\phi^n+\tau\beta a_{11}\gamma_1^{*}p^n+\tau\beta a_{12}\exp(\text{i}\mathcal L_h(c_1-c_2)\tau)\gamma_2^{*}p^n,\\
&\gamma_2^n=\exp(\text{i}\mathcal L_h c_2\tau)\phi^n+\tau\beta a_{21}\exp(\text{i}\mathcal L_h(c_2-c_1)\tau)\gamma_1^{*}p^n+\tau\beta a_{22}\gamma_2^{*}p^n.
\end{align*}
Multiply both sides of \eqref{NLS-fast-solver5} and \eqref{NLS-fast-solver6} with $\mathcal L_h$ and we then take discrete inner product with $\gamma_1^{*}$ and $\gamma_2^{*}$, respectively, to obtain
\begin{align}\label{NLS-fast-solver7}
&A_{11}\text{Re}\Big\langle\mathcal L_h\phi_{n1},\gamma_1^{*}\Big\rangle_h+A_{12}\text{Re}\Big\langle\mathcal L_h\phi_{n2},\gamma_2^{*}\Big\rangle_h=\text{Re}\Big\langle\mathcal L_h\gamma_1^n,\gamma_1^{*}\Big\rangle_h,\\\label{NLS-fast-solver8}
&A_{21}\text{Re}\Big\langle\mathcal L_h\phi_{n1},\gamma_1^{*}\Big\rangle_h+A_{22}\text{Re}\Big\langle\mathcal L_h\phi_{n2},\gamma_2^{*}\Big\rangle_h=\text{Re}\Big\langle\mathcal L_h\gamma_2^n,\gamma_2^{*}\Big\rangle_h,
\end{align} where
\begin{align}\label{NLS-fast-solver9}
&A_{11}=1+2\beta\tau^2\text{Re}\Big\langle a_{11}^2\mathcal L_h\gamma_1^{*}+a_{12}a_{21}\exp(\text{i}\mathcal L_h(c_1-c_2)\tau)\mathcal L_h\gamma_2^{*},\gamma_1^{*}\Big\rangle_h,\\\label{NLS-fast-solver10}
&A_{12}=2\beta\tau^2\text{Re}\Big\langle a_{11}a_{12}\mathcal L_h\gamma_1^{*}+a_{12}a_{22}\exp(\text{i}\mathcal L_h(c_1-c_2)\tau)\mathcal L_h\gamma_2^{*},\gamma_1^{*}\Big\rangle_h,\\\label{NLS-fast-solver11}
&A_{21}=2\beta\tau^2\text{Re}\Big\langle a_{21}a_{11}\exp(\text{i}\mathcal L_h(c_2-c_1)\tau)\mathcal L_h\gamma_1^{*}+a_{22}a_{21}\mathcal L_h\gamma_2^{*},\gamma_2^{*}\Big\rangle_h,\\\label{NLS-fast-solver12}
&A_{22}=1+2\beta\tau^2\text{Re}\Big\langle a_{21}a_{12}\exp(\text{i}\mathcal L_h(c_2-c_1)\tau)\mathcal L_h\gamma_1^{*}+a_{22}^2\mathcal L_h\gamma_2^{*},\gamma_2^{*}\Big\rangle_h.
\end{align}
Eqs. \eqref{NLS-fast-solver7} and \eqref{NLS-fast-solver8} form a $2\times 2$ linear system for the unknowns
$$\Big[\text{Re}\Big\langle\mathcal L_h\phi_{n1},\gamma_1^{*}\Big\rangle_h,\text{Re}\Big\langle\mathcal L_h\phi_{n2},\gamma_2^{*}\Big\rangle_h\Big]^T.$$

Solving $\Big[\text{Re}\Big\langle\mathcal L_h\phi_{n1},\gamma_1^{*}\Big\rangle_h,\text{Re}\Big\langle\mathcal L_h\phi_{n2},\gamma_2^{*}\Big\rangle_h\Big]^T$ from the $2\times 2$ linear system \eqref{NLS-fast-solver7} and \eqref{NLS-fast-solver8}, and $l_i,\ k_i$ and $\phi_{ni},\ i=1,2$ are then updated from \eqref{NLS-fast-solver2}-\eqref{NLS-fast-solver6}, respectively. Subsequently, $\phi^{n+1}$ and $p^{n+1}$ are obtained by \eqref{LI-E-scheme2}.

To summarize, 3rd-LI-IF-RKM and 4th-LI-IF-RKM can be efficiently implemented, respectively, as follows:
\begin{table}[H]
\tabcolsep=6pt
\small
\renewcommand\arraystretch{1.12}
\raggedleft
\begin{tabular*}{\textwidth}[h]{@{\extracolsep{\fill}}l}\hline 
 {\bf Algorithm 1: }  {\bf Scheme \ref{NLS-shn-3.1}}  of order 3 (denote by 3rd-LI-IF-RKM)\\[2ex]\hline
 {\bf Set the RK coefficients:}\\[2ex] {$c_1=\frac{1}{2}-\frac{\sqrt{3}}{6},\ c_2=\frac{1}{2}+\frac{\sqrt{3}}{6}$, $a_{11}=\frac{1}{4},\ a_{12}=\frac{1}{4}- \frac{\sqrt{3}}{6},\ a_{21}=\frac{1}{4}+ \frac{\sqrt{3}}{6},\ a_{22}=\frac{1}{4},\ b_1=\frac{1}{2},\ b_2=\frac{1}{2}$.}\\[2ex]
  {\bf Initialization:} apply the fourth-order exponential integration \cite{CWJ2021aml} to compute $\phi_{01},\phi_{02},\phi^1$ and $p^{1}$ \\[2ex]
 {\bf For n=1,2,$\cdots$, do}\\[2ex]
 Compute the  extrapolation approximation of order 3\\
 $\phi_{n1}^{*} = (2\sqrt3-4)\phi^{n-1}+(7\sqrt3)-11)\phi_{(n-1)1}+(6-5\sqrt3)\phi_{(n-1)2}+(10-4\sqrt3)\phi^{n}$\\[2ex]
$\phi_{n2}^{*} = -(2\sqrt3+4)\phi^{n-1}+(6+5\sqrt3)\phi_{(n-1)1}-(7\sqrt3+11)\phi_{(n-1)2}+(10+4\sqrt3)\phi^n$\\[2ex]
 Compute $\gamma_1^{*}$ and $\gamma_2^{*}$ by \eqref{NLS-fast-solver1}\\[2ex]
 Compute $A_{11},A_{12},A_{21},A_{22}$ by \eqref{NLS-fast-solver9}-\eqref{NLS-fast-solver12}\\[2ex]
 Compute $\text{Re}\Big\langle\mathcal L_h\phi_{n1},\gamma_1^{*}\Big\rangle_h$ and $\text{Re}\Big\langle\mathcal L_h\phi_{n2},\gamma_2^{*}\Big\rangle_h$ from \eqref{NLS-fast-solver7} and \eqref{NLS-fast-solver8}\\[2ex]
 Compute $l_i,\ k_i$ and $\phi_{ni},\ i=1,2$ from \eqref{NLS-fast-solver2}-\eqref{NLS-fast-solver6}, respectively\\[2ex]
 Update $\phi^{n+1}$ and $p^{n+1}$ via \eqref{LI-E-scheme2}
\\[2ex]
{\bf end for}\\ [2ex]
\hline
\end{tabular*}
\end{table}
\begin{rmk}{ As is stated in Remark \ref{NLS-rmk-3.4}, for our third-order scheme, we need to obtain two internal stage values on $t_0+c_1\tau$ and $t_0+c_2\tau$ (i.e., $\phi_{01}$ and $\phi_{02}$) for the following interpolation where $c_1$ and $c_2$ are the Gaussian quadrature nodes. Thus the fourth-order scheme is employed to obtain the initialization values in {\bf Algorithm 1}. Similarly, the sixth-order scheme is used to obtain the initialization values $\phi_{01},\ \phi_{02}$ and $\phi_{03}$ for 4th-LI-IF-RKM (see {\bf Algorithm 2}).}

\end{rmk}

\begin{table}[H]
\tabcolsep=6pt
\small
\renewcommand\arraystretch{1.12}
\raggedleft
\begin{tabular*}{\textwidth}[h]{@{\extracolsep{\fill}}l}\hline 
 {\bf Algorithm 2: }  {\bf Scheme \ref{NLS-shn-3.1}}  of order 4 (denote by 4th-LI-IF-RKM) \\[2ex]\hline
{{\bf Set the RK coefficients:}} \\[2ex]
{$c_1=\frac{1}{2}-\frac{\sqrt{15}}{10},\ c_2=\frac{1}{2},\ c_3=\frac{1}{2}+\frac{\sqrt{15}}{10}$,}\\[2ex]
{$a_{11}=\frac{5}{36},\ a_{12}=\frac{2}{9}-\frac{\sqrt{15}}{15},\ a_{13}=\frac{5}{36}- \frac{\sqrt{15}}{30},\ a_{21}=\frac{5}{36}+ \frac{\sqrt{15}}{24},\ a_{22}=\frac{2}{9},\ a_{23}=\frac{5}{36}- \frac{\sqrt{15}}{24}$,}\\[2ex]
{$a_{31}=\frac{5}{36}+ \frac{\sqrt{15}}{30},\ a_{32}=\frac{2}{9}+\frac{\sqrt{15}}{15},\ a_{33}=\frac{5}{36},\ b_1=\frac{5}{18},\ b_2=\frac{4}{9},\ b_3=\frac{5}{18}$}.\\[2ex]
  {\bf Initialization:} apply the sixth-order  exponential integration \cite{CWJ2021aml} to compute $\phi_{01},\phi_{02},\phi_{03},\phi^1$ and $p^{1}$ \\[2ex]
 {\bf For n=1,2,$\cdots$, do}\\[2ex]
 Compute the  extrapolation approximation of order 4 \\ $\phi_{n1}^{*}=(6\sqrt{15}-26)\phi^{n-1}+(11-5\sqrt{15}/3)\phi_{(n-1)1}+(16\sqrt{15}/3-24)\phi_{(n-1)2}+(40-29\sqrt{15}/3)\phi_{(n-1)3}$\\[2ex]
 $\phi_{n2}^{*}=-17\phi^{n-1}+(35/2+5\sqrt{15}/2)\phi_{(n-2)1}-17\phi_{(n-1)2}+(35/2-5\sqrt{15}/2)\phi_{(n-1)3}$\\[2ex]
 $\phi_{n3}^{*}=(-26-6\sqrt{15})\phi^{n-1}+(40+29\sqrt{15}/3)\phi_{(n-1)1}+(-24-16\sqrt{15}/3)\phi_{(n-1)2}+(11+5\sqrt{15}/3)\phi_{(n-1)3}$\\[2ex]
 Compute $\gamma_i^{*}=\text{i}\frac{|\phi_{ni}^{*}|^2\cdot\phi_{ni}^{*}}{\sqrt{\langle(\phi_{ni}^{*})^2,(\phi_{ni}^{*})^2\rangle_h+C_0}},\ i=1,2,3$\\[2ex]
 Compute\\
  $A_{11}=1+2\tau^2\beta\text{Re}\Big\langle a_{11}^2\mathcal L_h \gamma_1^{*}+a_{12}a_{21}\exp(\text i\mathcal L_h(c_1-c_2)\tau)\mathcal L_h\gamma_2^{*}+a_{13}a_{31}\exp(\text i\mathcal L_h(c_1-c_3)\tau)\mathcal L_h\gamma_3^{*},\gamma_1^{*}\Big\rangle_h$\\[2ex]
  $A_{12}=2\tau^2\beta\text{Re}\Big\langle a_{11}a_{12}\mathcal L_h \gamma_1^{*}+a_{12}a_{22}\exp(\text i\mathcal L_h(c_1-c_2)\tau)\mathcal L_h\gamma_2^{*}+a_{13}a_{32}\exp(\text i\mathcal L_h(c_1-c_3)\tau)\mathcal L_h\gamma_3^{*},\gamma_1^{*}\Big\rangle_h$\\[2ex]
  $A_{13}=2\tau^2\beta\text{Re}\Big\langle a_{11}a_{13}\mathcal L_h \gamma_1^{*}+a_{12}a_{23}\exp(\text i\mathcal L_h(c_1-c_2)\tau)\mathcal L_h\gamma_2^{*}+a_{13}a_{33}\exp(\text i\mathcal L_h(c_1-c_3)\tau)\mathcal L_h\gamma_3^{*},\gamma_1^{*}\Big\rangle_h$\\[2ex]
  $A_{21}=2\tau^2\beta\text{Re}\Big\langle a_{21}a_{11}\exp(\text i\mathcal L_h(c_2-c_1)\tau)\mathcal L_h \gamma_1^{*}+a_{22}a_{21}\mathcal L_h\gamma_2^{*}+a_{23}a_{31}\exp(\text i\mathcal L_h(c_2-c_3)\tau)\mathcal L_h\gamma_3^{*},\gamma_2^{*}\Big\rangle_h$\\ [2ex]
  $A_{22}=1+2\tau^2\beta\text{Re}\Big\langle a_{21}a_{12}\exp(\text i\mathcal L_h(c_2-c_1)\tau)\mathcal L_h \gamma_1^{*}+a_{22}^2\mathcal L_h\gamma_2^{*}+a_{23}a_{32}\exp(\text i\mathcal L_h(c_2-c_3)\tau)\mathcal L_h\gamma_3^{*},\gamma_2^{*}\Big\rangle_h$\\ [2ex]
  $A_{23}=2\tau^2\beta\text{Re}\Big\langle a_{21}a_{13}\exp(\text i\mathcal L_h(c_2-c_1)\tau)\mathcal L_h \gamma_1^{*}+a_{22}a_{23}\mathcal L_h\gamma_2^{*}+a_{23}a_{33}\exp(\text i\mathcal L_h(c_2-c_3)\tau)\mathcal L_h\gamma_3^{*},\gamma_2^{*}\Big\rangle_h$\\ [2ex]
 $A_{31}=2\tau^2\beta\text{Re}\Big\langle a_{31}a_{11}\exp(\text i\mathcal L_h(c_3-c_1)\tau)\mathcal L_h \gamma_1^{*}+a_{32}a_{21}\exp(\text i\mathcal L_h(c_3-c_2)\tau)\mathcal L_h\gamma_2^{*}+a_{33}a_{31}\mathcal L_h\gamma_3^{*},\gamma_3^{*}\Big\rangle_h$\\ [2ex]
 $A_{32}=2\tau^2\beta\text{Re}\Big\langle a_{31}a_{12}\exp(\text i\mathcal L_h(c_3-c_1)\tau)\mathcal L_h \gamma_1^{*}+a_{32}a_{22}\exp(\text i\mathcal L_h(c_3-c_2)\tau)\mathcal L_h\gamma_2^{*}+a_{33}a_{32}\mathcal L_h\gamma_3^{*},\gamma_3^{*}\Big\rangle_h$\\ [2ex]
 $A_{33}=1+2\tau^2\beta\text{Re}\Big\langle a_{31}a_{13}\exp(\text i\mathcal L_h(c_3-c_1)\tau)\mathcal L_h \gamma_1^{*}+a_{32}a_{23}\exp(\text i\mathcal L_h(c_3-c_2)\tau)\mathcal L_h\gamma_2^{*}+a_{33}^2\mathcal L_h\gamma_3^{*},\gamma_3^{*}\Big\rangle_h$\\ [2ex]
 $\gamma_1^n=\exp(\text i \mathcal L_hc_1\tau)\phi^n+\tau \beta a_{11}\gamma_1^{*}p^n+\tau\beta a_{12}\exp(\text i \mathcal L_h(c_1-c_2)\tau)\gamma_2^{*}p^n+\tau\beta a_{13}\exp(\text i \mathcal L_h(c_1-c_3)\tau)\gamma_3^{*}p^n$\\ [2ex]
 $\gamma_2^n=\exp(\text i \mathcal L_hc_2\tau)\phi^n+\tau \beta a_{21}\exp(\text i \mathcal L_h(c_2-c_1)\tau)\gamma_1^{*}p^n+\tau\beta a_{22}\gamma_2^{*}p^n+\tau\beta a_{23}\exp(\text i \mathcal L_h(c_2-c_3)\tau)\gamma_3^{*}p^n$\\ [2ex]
 $\gamma_3^n=\exp(\text i \mathcal L_hc_3\tau)\phi^n+\tau \beta a_{31}\exp(\text i \mathcal L_h(c_3-c_1)\tau)\gamma_1^{*}p^n+\tau\beta a_{32}\exp(\text i \mathcal L_h(c_3-c_2)\tau)\gamma_2^{*}p^n+\tau\beta a_{33}\gamma_3^{*}p^n$\\ [2ex]
 $B_1=\text{Re}\Big\langle\mathcal L_h\gamma_1^n,\gamma_1^{*}\Big\rangle_h,\ B_2=\text{Re}\Big\langle\mathcal L_h\gamma_2^n,\gamma_2^{*}\Big\rangle_h,\ B_3=\text{Re}\Big\langle\mathcal L_h\gamma_3^n,\gamma_3^{*}\Big\rangle_h$\\ [2ex]
Compute $\text{Re}\Big\langle\mathcal L_h\phi_{n1},\gamma_1^{*}\Big\rangle_h$, $\text{Re}\Big\langle\mathcal L_h\phi_{n2},\gamma_2^{*}\Big\rangle_h$  and $\text{Re}\Big\langle\mathcal L_h\phi_{n3},\gamma_3^{*}\Big\rangle_h$  from the
linear equation
$AX=B$ \\with coefficients of $[A_{ij}],\ [B_i],\ i,j=1,2,3$.\\[2ex]
 Compute $l_i,\ k_i$ and $\phi_{ni},\ i=1,2,3$ from \eqref{LI-E-scheme1}\\[2ex]
 Update $\phi^{n+1}$ and $p^{n+1}$ via \eqref{LI-E-scheme2}\\ [2ex]
{\bf end for}\\ [2ex]
\hline
\end{tabular*}
\end{table}

{\subsubsection{ Linearly implicit integrating
factor Runge-Kutta method based on the prediction-correction method}}
{\bf Scheme \ref{NLS-shn-3.1}} may cause large numerical errors and lack of robustness because of the extrapolation (see Section \ref{NLS-section-4.1}). Inspired by the previous works in \cite{BWC2020arixv,GZW2020jcp,SX-CiCP-2018}, we combine the prediction-correction method and the integrating factor Runge-Kutta method to propose a new class of high-order linearly implicit prediction-correction integrating factor Runge-Kutta methods (LI-PC-IF-RKMs). Here, we should note that the construction of such new scheme is only deferent from {\bf Scheme \ref{NLS-shn-3.1}} in obtaining $\phi_{ni}^*$, where the prediction-correction method is substituted for the extrapolation. For brevity, we directly give the scheme, as follows:
\begin{shm}\label{NLS-shn-3.2}  Let $b_i,a_{ij}\ (i,j=1,\cdots,s)$ be real numbers and $c_i=\sum_{j=1}^sa_{ij}$. For given $\phi^n$, $p^n$, an s-stage linearly implicit prediction-correction integrating factor Runge-Kutta method is given by

   1. Prediction: set $k_i^{(0)}=\phi^n,\ i=1,2,\cdots,s$ and let $M>0$ be a given integer, for $l=0,1,\cdots,M-1$, we calculate $k_i^{(l+1)}$ using
\begin{align}\label{LI-PC-scheme1}
\left\{
\begin{aligned}
&\phi_{ni}^{(l+1)}=\exp({\rm i}\mathcal{L}_hc_i\tau) \phi^n+\tau\sum_{j=1}^{s}a_{ij}\exp({\rm i}\mathcal{L}_h(c_i-c_j)\tau) k_j^{(l+1)},\\
&k_i^{(l+1)}={\rm i}\ \frac{\beta|\phi_{ni}^{(l)}|^2\cdot\phi_{ni}^{(l)} p_{ni}^{(l+1)}}
{\sqrt{\big\langle(\phi_{ni}^{(l)})^2,(\phi_{ni}^{(l)})^2\big\rangle_h\!+\!C_0}},\\
&p_{ni}^{(l+1)}=p^n+\tau\sum_{j=1}^{s}a_{ij}l_j^{(l+1)},\\
& l_i^{(l+1)}\!=2{\rm Re}\!\Big\langle{\rm i}\mathcal{L}_h\phi_{ni}^{(l)}, \frac{|\phi_{ni}^{(l)}|^2\cdot\phi_{ni}^{(l)}}
{\sqrt{\big\langle(\phi_{ni}^{(l)})^2,(\phi_{ni}^{(l)})^2\big\rangle_h\!+\!C_0}}\Big\rangle_h.
\end{aligned}
\right.
\end{align}
Then, we set $\phi_{ni}^{*}=\phi_{ni}^{(M)},\ i=1,2,\cdots,s$.

2.  Correction: for the predicted value $\phi_{ni}^{*}$, we obtain the $(\phi^{n+1},q^{n+1})$ by
\begin{align}\label{LI-PC-scheme2}
\left\{
\begin{aligned}
&\phi_{ni}=\exp({\rm i}\mathcal{L}_hc_i\tau) \phi^n+\tau\sum_{j=1}^{s}a_{ij}\exp({\rm i}\mathcal{L}_h(c_i-c_j)\tau) k_j,\\
&p_{ni}=p^n+\tau\sum_{j=1}^{s}a_{ij}l_j,\ l_i\!=2{\rm Re}\!\Big\langle\!{\rm i}\mathcal{L}_h\phi_{ni}, \frac{|\phi_{ni}^{*}|^2\cdot\phi_{ni}^{*}}
{\sqrt{\big\langle(\phi_{ni}^{*})^2,(\phi_{ni}^{*})^2\big\rangle_h\!+\!C_0}}\Big\rangle_h,\\
&k_i\!=\!{\rm i}\ \frac{\beta|\phi_{ni}^{*}|^2\cdot\phi_{ni}^{*} p_{ni}}
{\sqrt{\big\langle(\phi_{ni}^{*})^2,(\phi_{ni}^{*})^2\big\rangle_h\!+\!C_0}},
\end{aligned}
\right.
\end{align}
where $i=1,2,\cdots,s$.
Then $(\phi^{n+1},p^{n+1})$ are updated by
\begin{align}\label{LI-PC-scheme3}
\left\{
\begin{aligned}
&\phi^{n+1}=\exp({\rm i}\mathcal{L}_h\tau) \phi^n+\tau\sum_{i=1}^{s}b_i \exp({\rm i}\mathcal{L}_h(1-c_i)\tau) k_i,\\
&p^{n+1}=p^n+\tau\sum_{i=1}^sb_il_i.
\end{aligned}
\right.
\end{align}
\end{shm}By the similar argument as Theorem \ref{th2.1}, we can obtain:
\begin{thm} \label{NLS-PC-3.3} If the RK coefficients of {\bf Scheme \ref{NLS-shn-3.2}} satisfy \eqref{LI-EI-coefficient}
 {\bf Scheme \ref{NLS-shn-3.2}} can preserve the discrete quadratic energy \eqref{LI-EI-F-energy}.
\end{thm}
{In addition, if we choose the RK coefficients of the $s$-stage Gauss collocation method and set a proper iterative steps $M$, {\bf Scheme \ref{NLS-shn-3.2}} can achieve $2s$-order accuracy. In particular, it can also be solved efficiently.} For simplicity, we take $s=2$ as an example, as follows:
\begin{table}[H]
\tabcolsep=6pt
\small
\renewcommand\arraystretch{1.12}
\raggedleft
\begin{tabular*}{\textwidth}[h]{@{\extracolsep{\fill}}l}\hline 
 {\bf Algorithm 3: }  {\bf Scheme \ref{NLS-shn-3.2}} with 2-stage Gauss collocation method \\[2ex]\hline
 {\bf Set the RK coefficients:}\\[2ex] {$c_1=\frac{1}{2}-\sqrt{3}{6},\ c_2=\frac{1}{2}+\frac{3}{6}$, $a_{11}=\frac{1}{4},\ a_{12}=\frac{1}{4}- \frac{\sqrt{3}}{6},\ a_{21}=\frac{1}{4}+ \frac{\sqrt{3}}{6},\ a_{22}=\frac{1}{4},\ b_1=\frac{1}{2},\ b_2=\frac{1}{2}$.}\\[2ex]
 {\bf For n=1,2,$\cdots$, do}\\[2ex]
 Set $M$, compute the $\phi_{n1}^{*},\ \phi_{n2}^{*}$ by using \eqref{LI-PC-scheme1}\\[2ex]
 Compute $\gamma_1^{*}$ and $\gamma_2^{*}$ by \eqref{NLS-fast-solver1}\\[2ex]
 Compute $A_{11},A_{12},A_{21},A_{22}$ by \eqref{NLS-fast-solver9}-\eqref{NLS-fast-solver12}\\[2ex]
 Compute $\text{Re}\Big\langle\mathcal L_h\phi_{n1},\gamma_1^{*}\Big\rangle_h$ and $\text{Re}\Big\langle\mathcal L_h\phi_{n2},\gamma_2^{*}\Big\rangle_h$ from \eqref{NLS-fast-solver7} and \eqref{NLS-fast-solver8}\\[2ex]
 Compute $l_i,\ k_i$ and $\phi_{ni},\ i=1,2$ from \eqref{NLS-fast-solver2}-\eqref{NLS-fast-solver6}, respectively\\[2ex]
 Update $\phi^{n+1}$ and $p^{n+1}$ via \eqref{LI-E-scheme2}
\\[2ex]
{\bf end for}\\ [2ex]
\hline
\end{tabular*}
\end{table}
 \begin{rmk}  {To the best of our knowledge,
the theoretical result on the choice of the iterative steps $M$ is missing. From our numerical
experiences, if we take the iterative steps $M=4$ and the RK coefficients of 2-stage Gauss collocation method, {\bf Scheme \ref{NLS-shn-3.2}} can be fourth-order accuracy in time and  the resulting scheme is denote by 4th-LI-PC-IF-RKM.}
\end{rmk}

\begin{rmk} \label{LI-EI-rmk3.4} We consider a Hamiltonian PDEs system given by
\begin{align}\label{Hamiltionian-orginal-system}
\partial_t z = \mathcal{D}(\mathcal A z+F^{'}(z)),\ z\in\Omega\subset\mathbb R^d,\ t> 0, d=1,2,3,
\end{align}
with a Hamiltonian energy
 \begin{align}\label{H-energy:2.1}
 \mathcal{E}(t)=\frac{1}{2}(z,\mathcal Az) + \big(F(z),1\big),\ t\ge0,
 \end{align}
where $z := z({\bf x},t),$ $\mathcal{D}$ is a constant skew-adjoint operation and $\mathcal A$ is a self-adjoint operation and $\int_{\Omega}F(z({\bf x},t))d{\bf x}$ is bounded from below.
By introducing a scalar auxiliary variable
\begin{align*}
p:=p(t)=\sqrt{2(F(z),1)+2C_0},
\end{align*}
 where $C_0$ is a constant large enough to make $p$ well-defined for all $z({\bf x},t),\ {\bf x}\in\Omega$. On the basis of the energy-variational principle, the system \eqref{Hamiltionian-orginal-system}  can be reformulated into
 \begin{align}\label{IEQ-Formulation}
 \left\lbrace
  \begin{aligned}
&\partial_t z = \mathcal{D}\Big(\mathcal Az + \frac{F^{'}(z)}{\sqrt{2\Big(F(z),1\Big)+2C_0}} p\Big),\\
&\partial_t p = \frac{\Big(F^{'}(z),\partial_t z\Big)}{\sqrt{2\Big(F(z),1\Big)+2C_0}},
 \end{aligned}\right.
  \end{align}
  which preserves the following  quadratic energy
  \begin{align}\label{Q-energy-conservetion-law}
  \frac{d}{dt} \mathcal{E}(t)=0,\ \mathcal{E}(t)=\frac{1}{2}(z,\mathcal Az)+\frac{1}{2}p^2-C_0,\ t\ge 0.
  \end{align}
  The proposed linearly implicit integrating factor Runge-Kutta methods
can be easily generalized to solve the above system.
\end{rmk}

\begin{rmk}When the standard Fourier pseudo-spectral method is employed to the system \eqref{NLS-equation} for spatial discretizations, we can obtain
discrete Hamiltonian energy at $t=t_n$ as
\begin{align}\label{sm-Hamiltonian-energy-conservetion-law}
\mathcal E_h^n =\langle\phi^n,\mathcal L_h\phi^n\rangle_h+\frac{\beta}{2}h_xh_y\sum_{j=0}^{N_x-1}\sum_{k=0}^{N_y-1}|\phi_{j,k}^n|^4,\ n=0,1,2,\cdots,.
\end{align}
We note that however, the modified energy introduced by the SAV approach is only equivalent to the Hamiltonian energy in the continuous sense, but not for the discrete sense. Thus, the proposed scheme cannot preserve discrete Hamiltonian energy \eqref{sm-Hamiltonian-energy-conservetion-law}. In addition, the nonlinear terms of the semi-discrete system \eqref{NLS-SP-SAV-reformulation} are explicitly discretized, thus the following discrete mass also cannot be preserved exactly by the proposed scheme
\begin{align}\label{LI-EI-discerete-mass}
\mathcal M_h^n=h_xh_y\sum_{j=0}^{N_x-1}\sum_{k=0}^{N_y-1}|\phi_{j,k}^n|^2,\ n=0,1,2,\cdots,.
\end{align}

\end{rmk}

\begin{rmk} { In Refs. \cite{JLQYarix2020,JLQZmc18}, Ju et al. presented the convergence analysis of a second-order energy stable ETD scheme for the epitaxial growth model without slope
selection. More recently, Akrivis et al. \cite{ALL2019sisc} established the optimal error estimate of a class of linear high-order energy stable RK methods for the Allen-Cahn and Cahn-Hilliard equations. Thus, how to combine their analysis techniques to analyze our schemes will be
considered in the future works.}

\end{rmk}

\section{Numerical results}\label{Sec:LI-EP:4}
In this section, we report the numerical performances in the accuracy, computational efficiency and
invariants-preservation of the proposed schemes of the NLSE. For brevity, in the rest
of this paper,  3rd-LI-IF-RKM  (see {\bf Algorithm 1}), 4th-LI-IF-RKM (see {\bf Algorithm 2}) and  4th-LI-PC-IF-RKM (see {\bf Algorithm 3}) are
only used for demonstration purposes. In order to quantify the numerical solution, we use the $L^2$- and $L^{\infty}$-norms of the error between the numerical solution $\phi_{j,k}^n$ and the exact solution
$\phi(x_j,y_k,t_n)$, respectively. {Meanwhile, we also investigate the relative residuals on the mass, Hamiltonian energy and modified energy, defined respectively, as
\begin{align*}
RM^n=\Big|\frac{\mathcal M_h^n-\mathcal M_h^0}{\mathcal M_h^0}\Big|,\ RH^n=\Big|\frac{\mathcal E_h^n-\mathcal E_h^0}{\mathcal E_h^0}\Big|,\ RE^n=\Big|\frac{E_h^n-E_h^0}{ E_h^0}\Big|.
\end{align*}}
Furthermore, we compare the selected schemes with some existing structure-preserving schemes, as follows:
{\begin{itemize}
\item SAVS: classical SAV scheme in \cite{AST2019P};
\item M-CNS: modified Crank-Nicolson scheme in \cite{FLM2019cicp};
\item M-BDF2: modified BDF scheme ($k=2$) in \cite{FLM2019cicp};
\item 3rd-LI-EPS: linearly implicit energy-preserving scheme of order 3 in \cite{LGZ2020aml};
\item 4th-Gauss: Gauss method of order 4 in \cite{FLM2021SIUM};
\item 4th-LI-EPS: linearly implicit energy-preserving scheme of order 4 in \cite{LGZ2020aml};
\end{itemize}
As a summary, the properties of these schemes are displayed in Tab. \ref{Tab-IPs:1}.
\begin{table}[H]
\tabcolsep=6pt
\small
\renewcommand\arraystretch{1.2}
\centering
\caption{Comparisons of properties of different numerical schemes}\label{Tab-IPs:1}
\begin{tabular*}{\textwidth}[h]{@{\extracolsep{\fill}}c c c c c c c c c c c c c  }\hline 
 \diagbox{Scheme}{Property}& SAVS& M-CNS& M-BDF2 &3rd-LI-EPS&  3rd-LI-IF-RKM\\\hline
  Conserving energy \eqref{sm-Hamiltonian-energy-conservetion-law}& No&Yes &No& No &No\\[1ex]
  Conserving mass \eqref{LI-EI-discerete-mass}& No&Yes &No& Yes &No\\[1ex]
 Conserving energy \eqref{LI-EI-F-energy}&Yes&No &No& No &Yes\\[1ex]
 Temporal accuracy& 2nd&2nd&2nd&3rd&3rd \\[1ex]
 Fully implicit &No&Yes&Yes&No&No\\
 Linearly implicit &Yes&No&No&Yes&Yes\\
\hline
\end{tabular*}
\vspace{0.5mm}
\centering
\caption*{}
\begin{tabular*}{\textwidth}[h]{@{\extracolsep{\fill}}c c c c c c c c c c c c c  }\hline 
 \diagbox{Scheme}{Property}&4th-Gauss& 4th-LI-EPS&4th-LI-PC-IF-RKM& 4th-LI-IF-RKM \\\hline
  Conserving energy \eqref{sm-Hamiltonian-energy-conservetion-law}& No& No &No&No \\[1ex]
  Conserving mass \eqref{LI-EI-discerete-mass}& Yes& Yes &No&No \\[1ex]
 Conserving energy \eqref{LI-EI-F-energy}& Yes& No &Yes&Yes\\[1ex]
 Temporal accuracy& 4th&4th&4th&4th \\[1ex]
 Fully implicit &Yes&No&No&No\\
 Linearly implicit &No&Yes&Yes&Yes\\
\hline
\end{tabular*}

\end{table}

We should note that the discretizations of all the selected schemes in space are the standard Fourier pseudo-spectral method.
In our computations, the standard fixed-point iteration is used for the fully implicit schemes and the Jacobi iteration method is employed for the linear system obtained by 3rd-LI-EPS and 4th-LI-EPS. Here, the iteration will terminate if the infinity norm of the error between two adjacent iterative steps is less than $10^{-14}$. Furthermore, the FFT technique is employed to
compute efficiently the exponential matrices and the constant coefficient matrices at every nonlinear iteration step.}

\subsection{1D nonlinear Schr\"odinger equation}\label{NLS-section-4.1}
We first use the proposed schemes to solve the one dimensional (1D) nonlinear Schr\"{o}dinger equation \eqref{NLS-equation} with $\beta=2$, which admits a single soliton solution given by \cite{ABB13}
\begin{align}
\phi(x,t)=e^{\text{i}(2x-3t)}\text{sech}(x-4t),\ x\in\Omega\subset\mathbb{R},\ t\ge0.
\end{align}

The computations are done on the space interval $\Omega=[-40,40]$ and  a periodic boundary is considered. To test the temporal discretization errors of the numerical schemes, we fix the Fourier node $1024$ so that spatial errors play no role here.

{ The $L^2$ errors and $L^{\infty}$ errors in numerical solution of $\phi$ at $T=1$ are calculated by using nine numerical schemes with various time steps $\tau=\frac{0.01}{2^i},\ i=0,1,2,4$ and the results are summarized in Fig. \ref{LI-EI-1d-scheme:fig:1}, which shows SAVS, M-CNS and M-BDF2 have second-order accuracy in time, and 3rd-LI-IF-RKM and 3rd-LI-EPS have third-order accuracy in time, and 4th-Gauss, 3rd-LI-IF-RKM, 4th-LI-PC-IF-RKM and 4th-LI-EPS have fourth-order accuracy in time, respectively. In addition, at the same time steps, we can observe that (i) the $L^2$ errors and $L^{\infty}$ errors of the high-order schemes are significantly smaller than the second-order schemes; (ii) 4th-LI-PC-IF-RKM provides the smallest errors, and the errors provided by 3rd-LI-IF-RKM are much smaller than the ones provided by 3rd-LI-EPS; (iii) the errors provided by 4th-Gauss are much smaller than the ones provided by 4th-LI-IF-RKM and we thank that the large numerical error is introduced by the extrapolation approximation.

\begin{figure}[H]
\centering\begin{minipage}[t]{60mm}
\includegraphics[width=60mm]{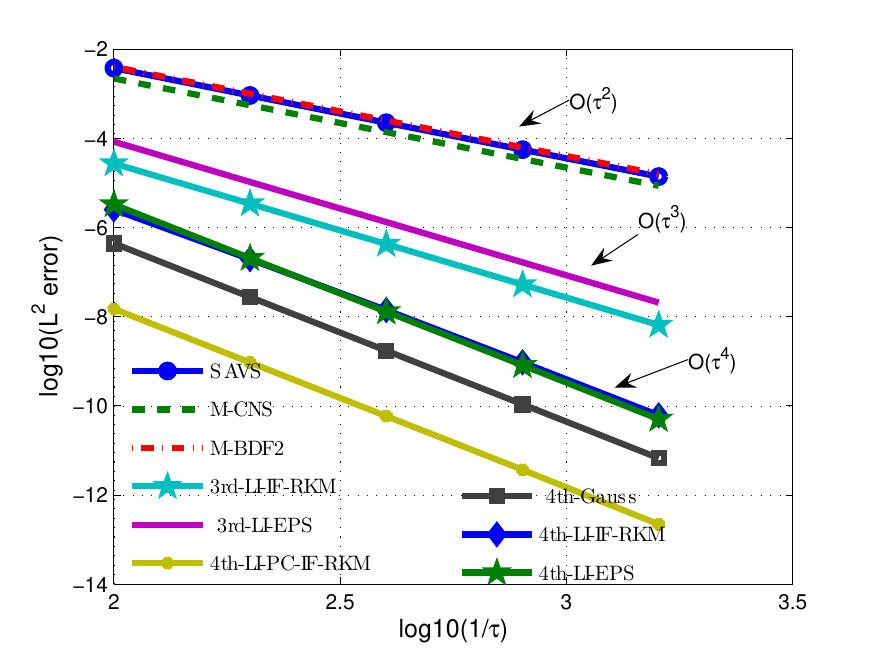}
\end{minipage}
\begin{minipage}[t]{60mm}
\includegraphics[width=60mm]{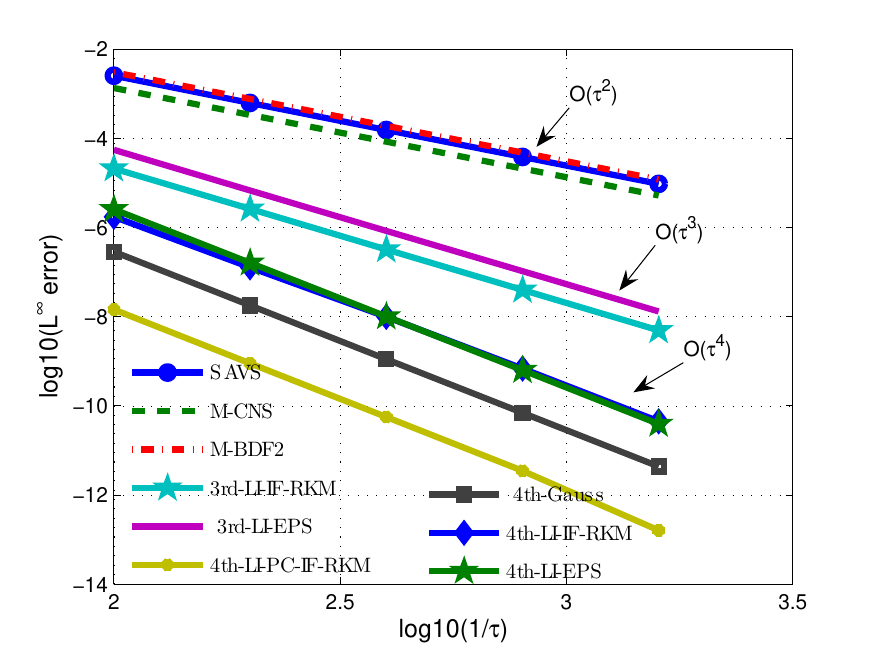}
\end{minipage}
\caption{ Time step refinement tests using the different numerical schemes for the 1D Schr\"odinger equation \eqref{NLS-equation}.}\label{LI-EI-1d-scheme:fig:1}
\end{figure}

Then, to demonstrate the advantages of the proposed high-order schemes with the existing schemes, we investigate the global $L^2$ errors and $L^{\infty}$ errors of $\phi$ versus the CPU costs using the different schemes at $T=1$ with the Fourier node 2048. The results are summarized in Fig. \ref{LI-EI-1d-scheme:fig:2}. As it is seen, for a given global error, we can draw the following observations: (i) the high-order schemes spend much less CPU time than all of second-order schemes and the cost of 4th-LI-PC-IF-RKM is cheapest; (ii) the costs of 3rd-LI-IF-RKM and 4th-LI-IF-RKM are much cheaper than 3rd-LI-EPS and 4th-LI-EPS, respectively; (iii) the cost of 4th-LI-IF-RKM is a little expensive over 4th-Gauss. We thank that on the one hand, at the same time steps, the errors provided by 4th-Gauss are much smaller than 4th-LI-IF-RKM (see Figure \ref{LI-EI-1d-scheme:fig:1}). On the other hand, at every nonlinear iteration step, the constant coefficient matrix of 4th-Gauss is explicitly solved based on the matrix diagonalization method \cite{ST06}  and FFT, thus it is very cheap.

\begin{figure}[H]
\centering\begin{minipage}[t]{60mm}
\includegraphics[width=60mm]{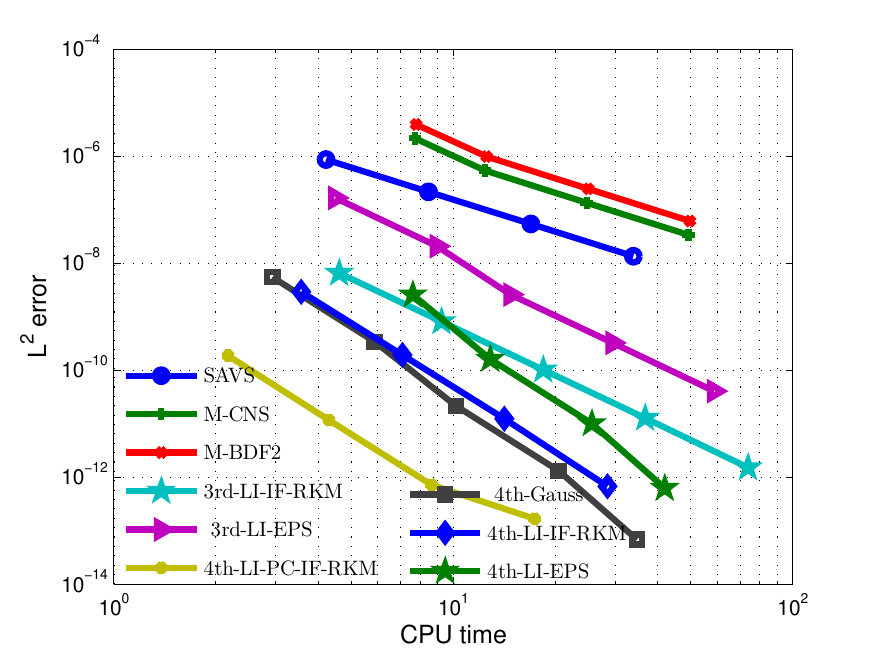}
\end{minipage}
\begin{minipage}[t]{60mm}
\includegraphics[width=60mm]{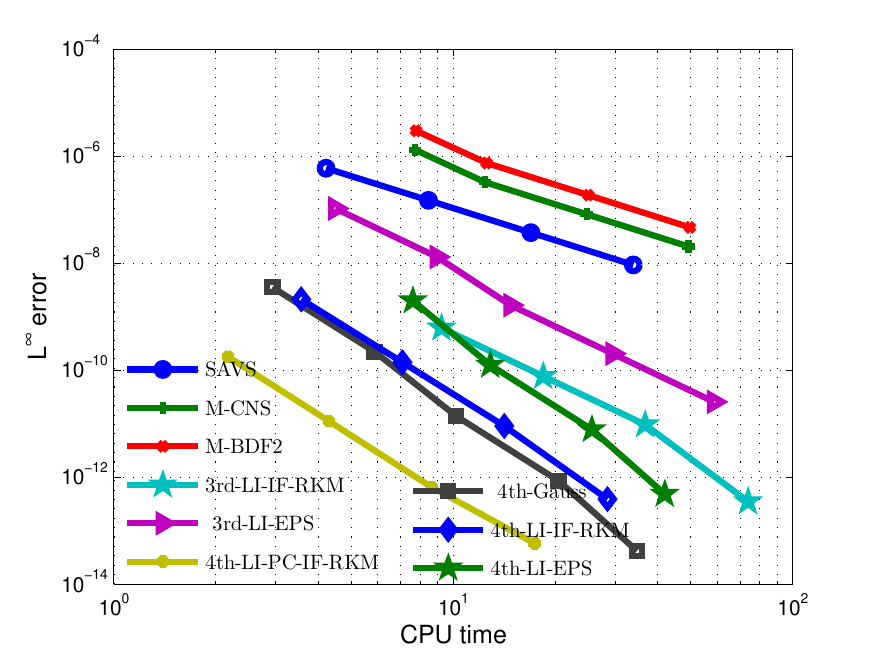}
\end{minipage}
\caption{ The numerical errors versus the CPU costs using the different numerical schemes for the 1D Schr\"odinger equation \eqref{NLS-equation}.}\label{LI-EI-1d-scheme:fig:2}
\end{figure}

Subsequently, we investigate the robustness of the proposed schemes in simulating evolution of the soliton as a large time step $\tau=0.1$ is chosen.
The results are summarized in Figs. \ref{LI-EI-1d-scheme:fig:3}-\ref{LI-EI-1d-scheme:fig:4}, where the result provided by 4th-LI-EPS is omitted since the computing result is below-up. It is clear to see that (i)  the computational results provided by SAVS, 3rd-LI-IF-RKM and 4th-LI-IF-RKM are wrong;  (ii) although M-BDF2, M-CNS and 3rd-LI-EPS can capture the amplitude and
waveform, but the small disturbance emergences; (iii) 4th-LI-PC-IF-RKM and 4th-Gauss simulate the soliton well, however, when we set the time step as $\tau=0.2$, the small disturbance emergences for 4th-Gauss (see Fig. \ref{LI-EI-1d-scheme:fig:4} (a)), which implies that 4th-LI-PC-IF-RKM is more robust. In addition, based on the results of 3rd-LI-EPS, 4th-LI-EPS, 3rd-LI-IF-RKM and 4th-LI-IF-RKM, we should note that the instability introduced by the extrapolation approximation will be enhanced as the degree of the interpolation polynomial increases, and scuh instability will probably make the proposed scheme invalid as the large time step is chosen. Thus, how to construct an accurate and stable extrapolation approximation need to be further discussed.

\begin{figure}[H]
\begin{minipage}[t]{60mm}
\includegraphics[width=60mm]{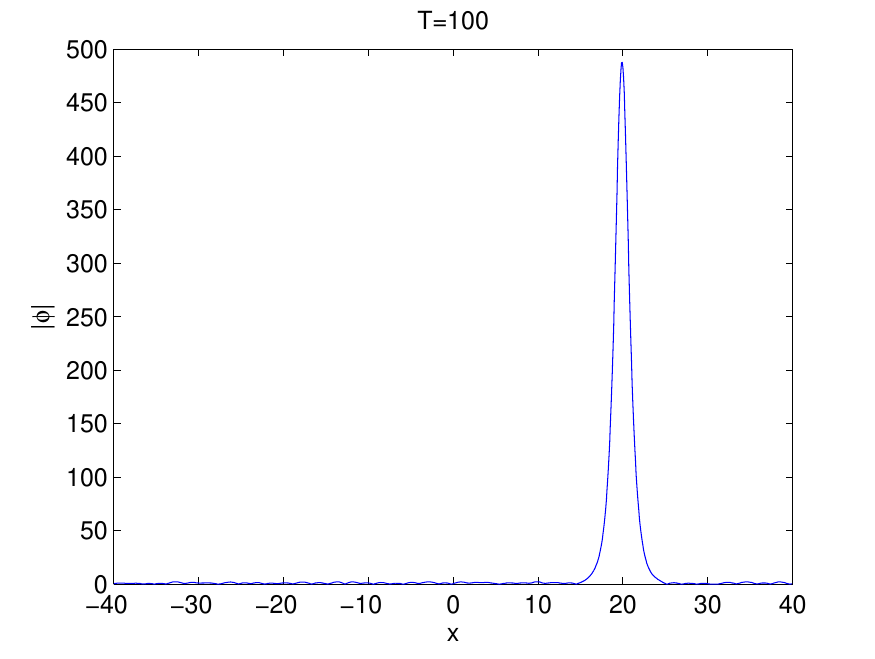}
\caption*{(a) SAVS }
\end{minipage}
\begin{minipage}[t]{60mm}
\includegraphics[width=60mm]{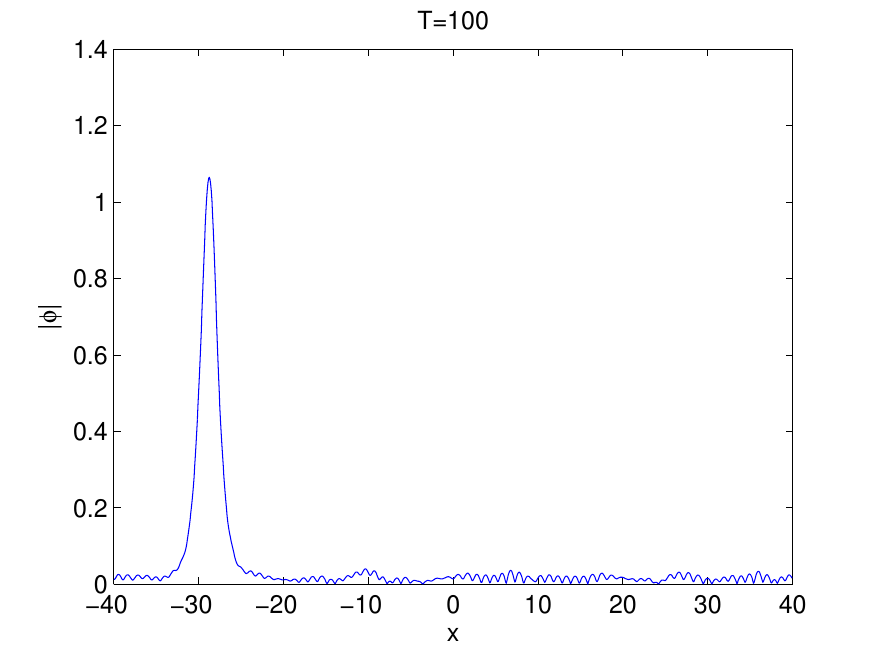}
\caption*{(b)  M-BDF2}
\end{minipage}
\centering\begin{minipage}[t]{60mm}
\includegraphics[width=60mm]{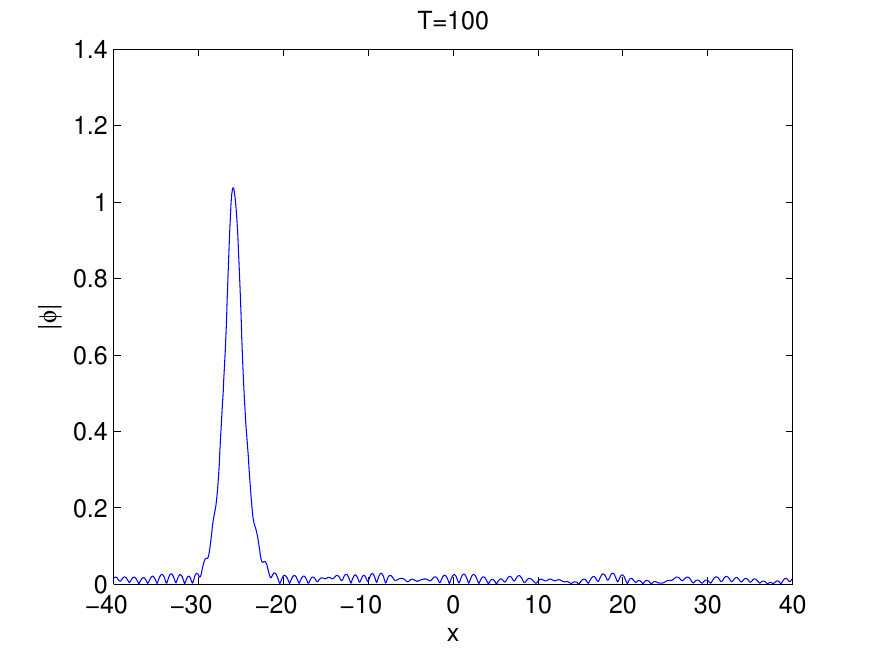}
\caption*{(c)  M-CNS}
\end{minipage}
\begin{minipage}[t]{60mm}
\includegraphics[width=60mm]{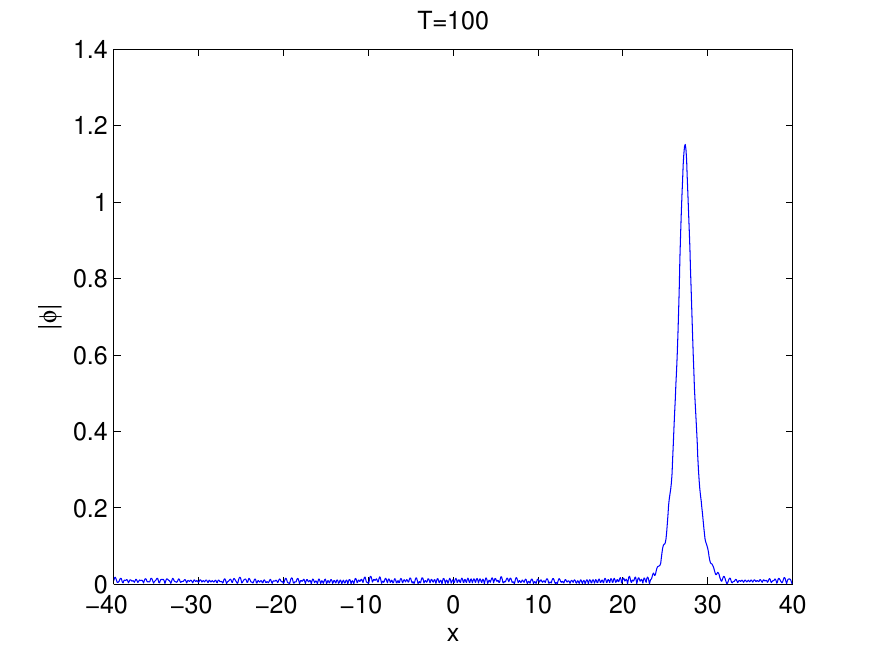}
\caption*{(d)  3rd-LI-EPS}
\end{minipage}
\begin{minipage}[t]{60mm}
\includegraphics[width=60mm]{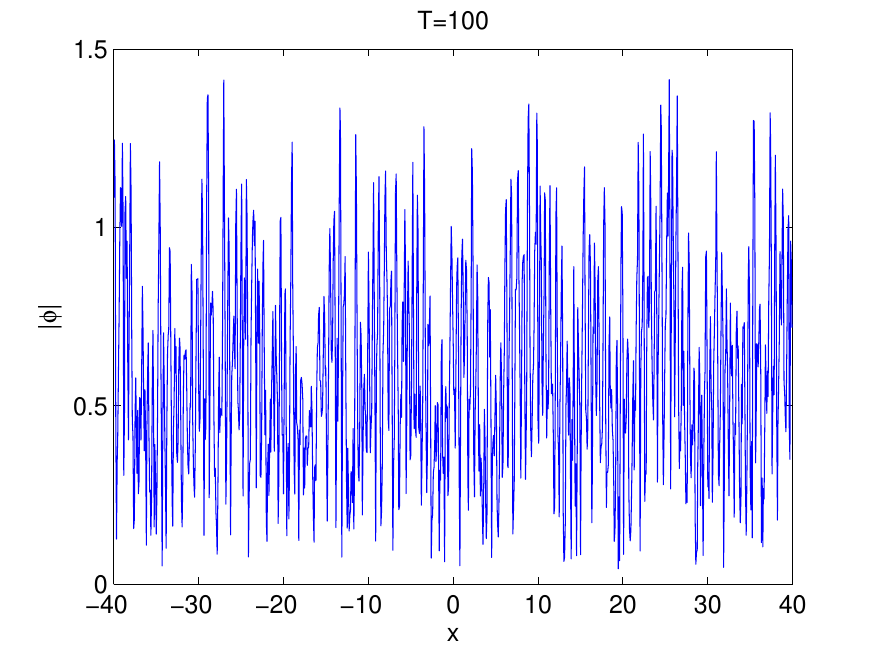}
\caption*{(e)  3rd-LI-IF-RKM}
\end{minipage}
\begin{minipage}[t]{60mm}
\includegraphics[width=60mm]{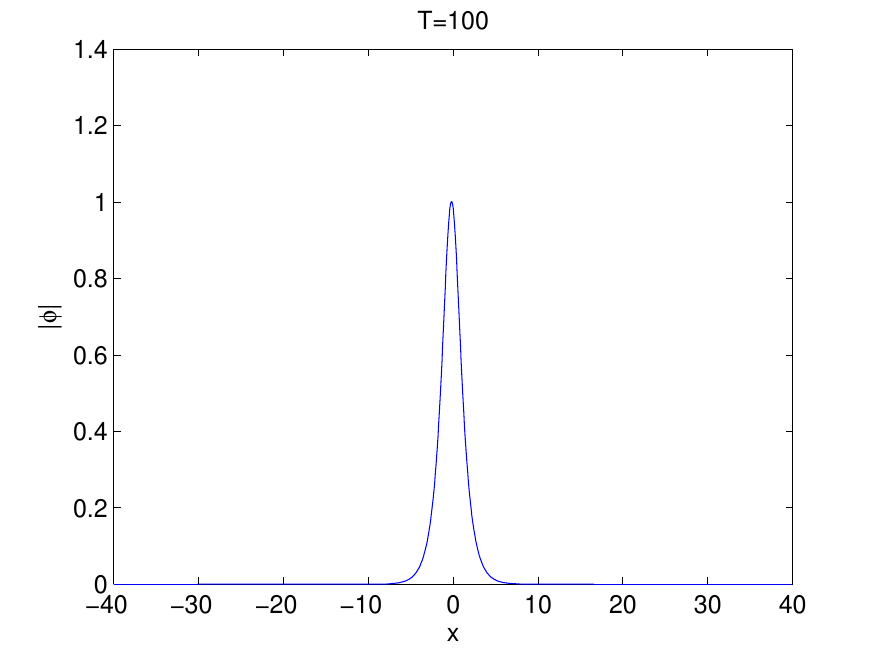}
\caption*{(f)  4th-Gauss}
\end{minipage}
\end{figure}
\begin{figure}[H]
\centering\begin{minipage}[t]{60mm}
\includegraphics[width=60mm]{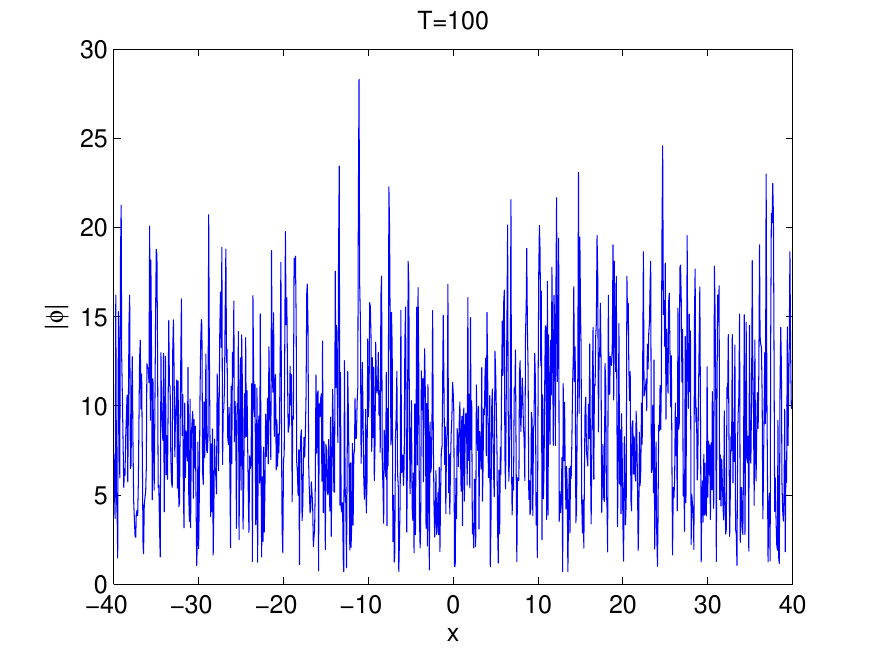}
\caption*{(g)  4th-LI-IF-RKM}
\end{minipage}
\begin{minipage}[t]{60mm}
\includegraphics[width=60mm]{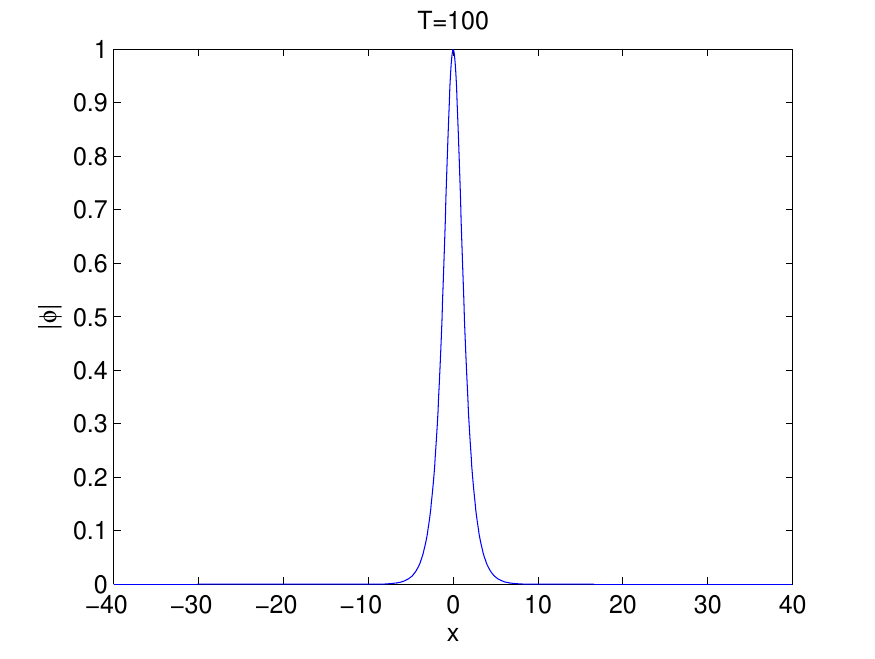}
\caption*{(h)   4th-LI-PC-IF-RKM}
\end{minipage}
\caption{ The numerical solution at $T=100$ using the different numerical schemes for the 1D Schr\"odinger equation \eqref{NLS-equation}. We take the parameter $\beta=2$, the time step $\tau=0.1$ and Fourier node $N=1024$ and the initial condition is chosen as the exact solution at $t=0$.}\label{LI-EI-1d-scheme:fig:3}
\end{figure}
\begin{figure}[H]
\centering\begin{minipage}[t]{60mm}
\includegraphics[width=60mm]{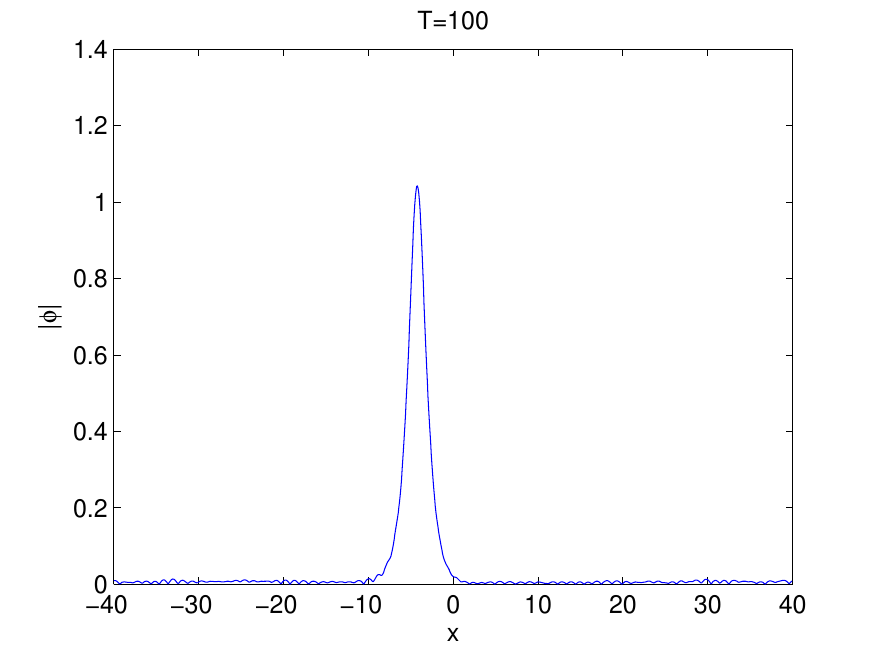}
\caption*{(a) 4th-Gauss}
\end{minipage}
\begin{minipage}[t]{60mm}
\includegraphics[width=60mm]{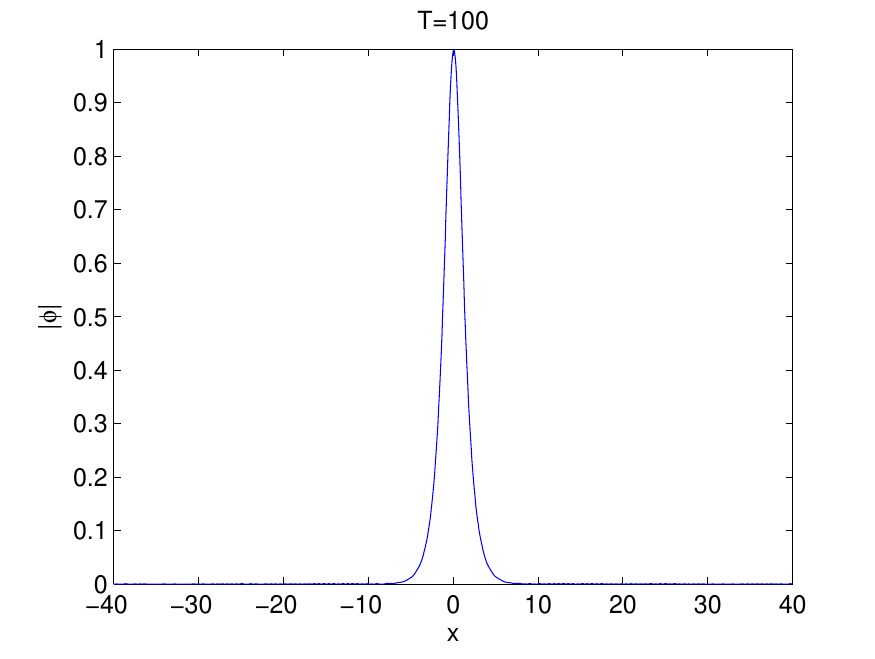}
\caption*{(b)  4th-LI-PC-IF-RKM}
\end{minipage}
\caption{ The numerical solution at $T=100$ using the two different numerical schemes for the 1D Schr\"odinger equation \eqref{NLS-equation}. We take the parameter $\beta=2$, the time step $\tau=0.2$ and Fourier node $N=1024$ and the initial condition is chosen as the exact solution at $t=0$.}\label{LI-EI-1d-scheme:fig:4}
\end{figure}}

\begin{figure}[H]
\centering
\begin{minipage}[t]{45mm}
\includegraphics[width=45mm]{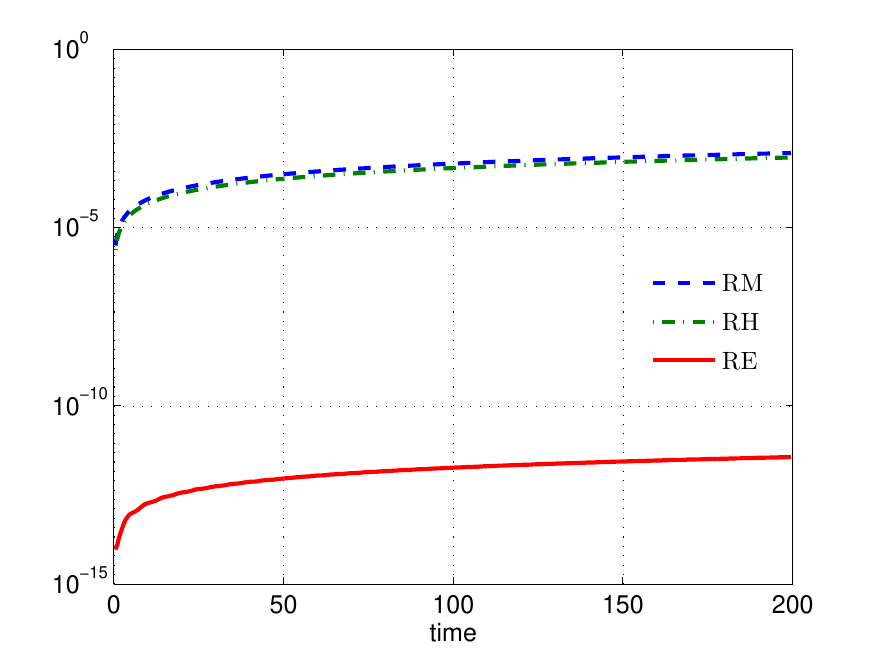}
\caption*{(a)  3rd-LI-IF-RKM}
\end{minipage}
\begin{minipage}[t]{45mm}
\includegraphics[width=45mm]{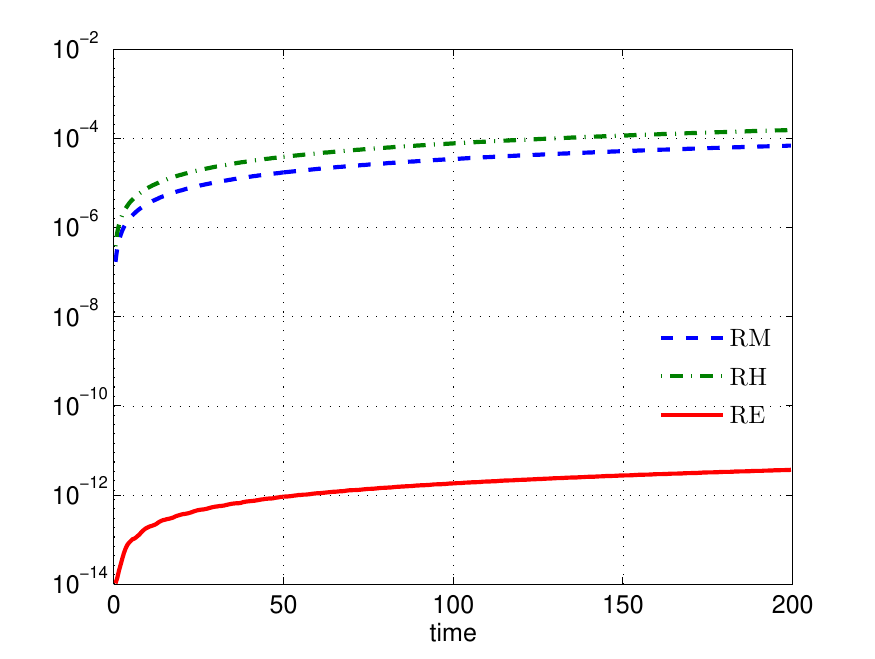}
\caption*{(b)  4th-LI-IF-RKM}
\end{minipage}
\begin{minipage}[t]{45mm}
\includegraphics[width=45mm]{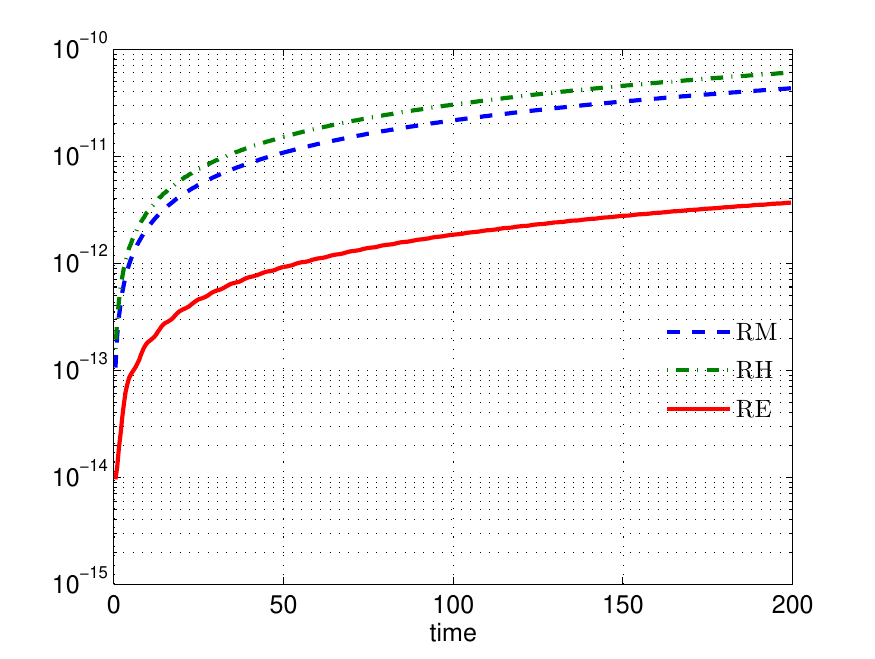}
\caption*{(c)  4th-LI-PC-IF-RKM}
\end{minipage}
\caption{The relative residuals on three discrete conservation laws using the proposed schemes for the 1D Schr\"odinger equation \eqref{NLS-equation}. We take the parameter $\beta=2$, time step $\tau=0.01$ and Fourier node 256 and the initial condition is chosen as the exact solution at $t=0$.}\label{LI-EI-1d-scheme:fig:5}
\end{figure}

Finally, we visualize the relative residuals on the Hamiltonian energy \eqref{sm-Hamiltonian-energy-conservetion-law}, mass \eqref{LI-EI-discerete-mass}
 and quadratic energy \eqref{LI-EI-F-energy}, which are summarized in Fig. \ref{LI-EI-1d-scheme:fig:5}. We can observe that 3rd-LI-IF-RKM, 4th-LI-IF-RKM and 4th-LI-PC-IF-RKM  preserve the quadratic energy exactly and the discrete mass and Hamiltonian energy can be preserved well.

\subsection{2D nonlinear Schr\"odinger equation}
In this subsection, we employ the proposed schemes to solve the two dimensional (2D) nonlinear Schr\"odinger equation, which possesses the following analytical solution
\begin{align*}
\phi(x,y,t)=A\exp(\text{i}(k_1x+k_2y-\omega t)),\ \omega=k_1^2+k_2^2-\beta|A|^2,\ (x,y)\in\Omega.
\end{align*}

\begin{figure}[H]
\centering\begin{minipage}[t]{60mm}
\includegraphics[width=60mm]{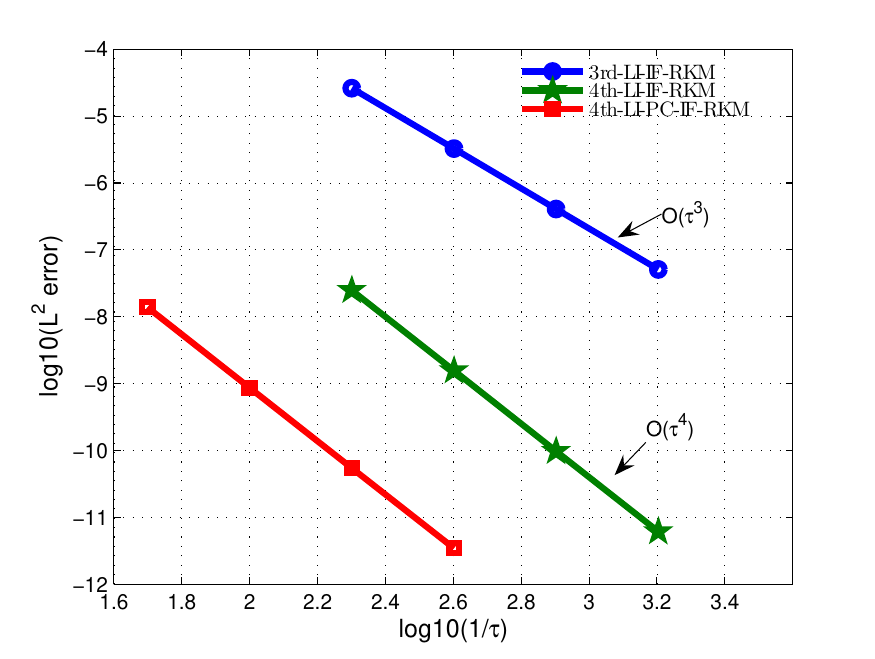}
\end{minipage}
\begin{minipage}[t]{60mm}
\includegraphics[width=60mm]{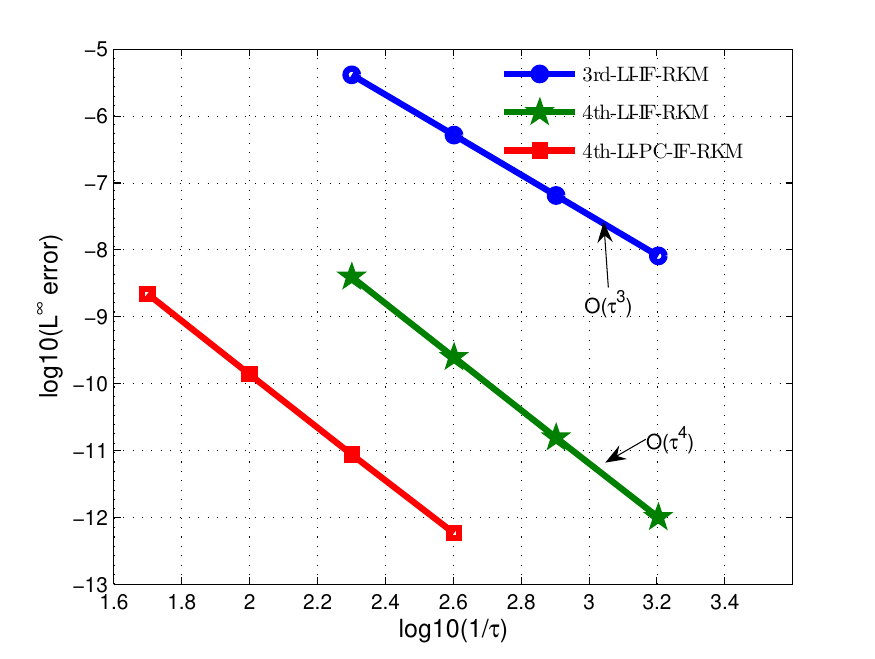}
\end{minipage}
\caption{ Time step refinement tests using the three proposed numerical methods for the 2D Schr\"odinger equation \eqref{NLS-equation}.}\label{LI-EI-2d-scheme:fig:1}
\end{figure}

\begin{figure}[H]
\begin{minipage}[t]{45mm}
\includegraphics[width=45mm]{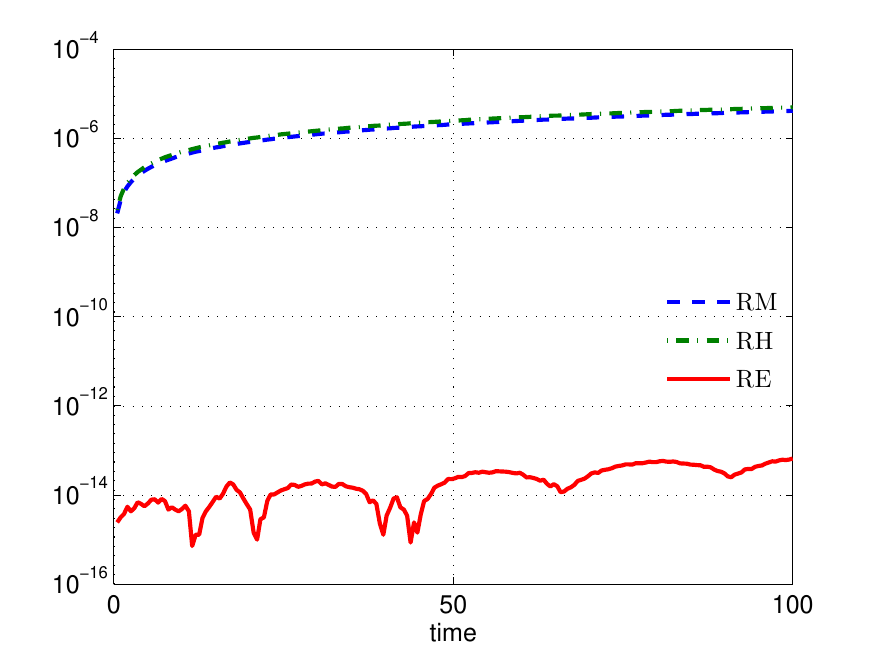}
  \caption*{(a)  3rd-LI-IF-RKM}
\end{minipage}
\begin{minipage}[t]{45mm}
\includegraphics[width=45mm]{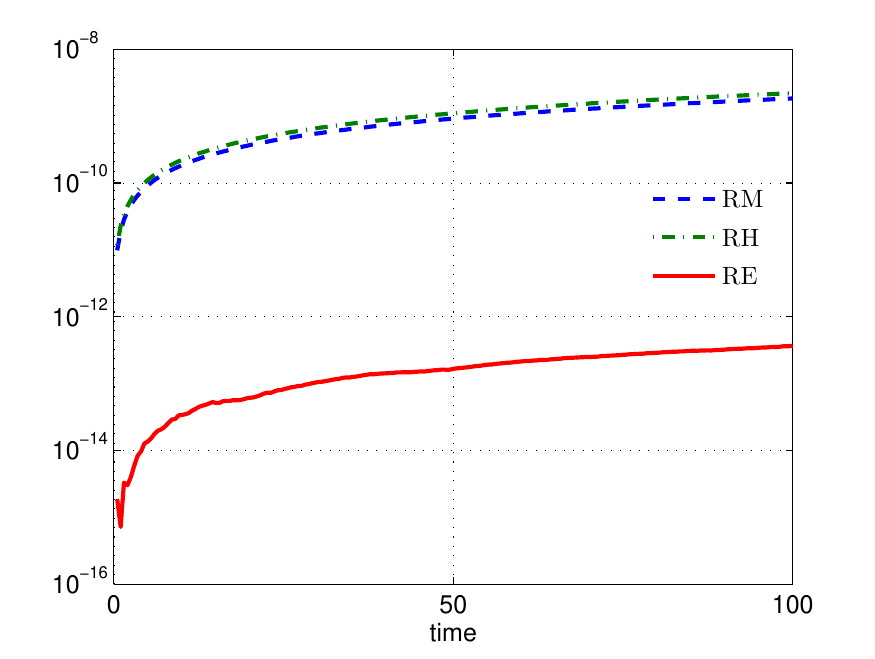}
  \caption*{(b)  4th-LI-IF-RKM}
\end{minipage}
\centering\begin{minipage}[t]{45mm}
\includegraphics[width=45mm]{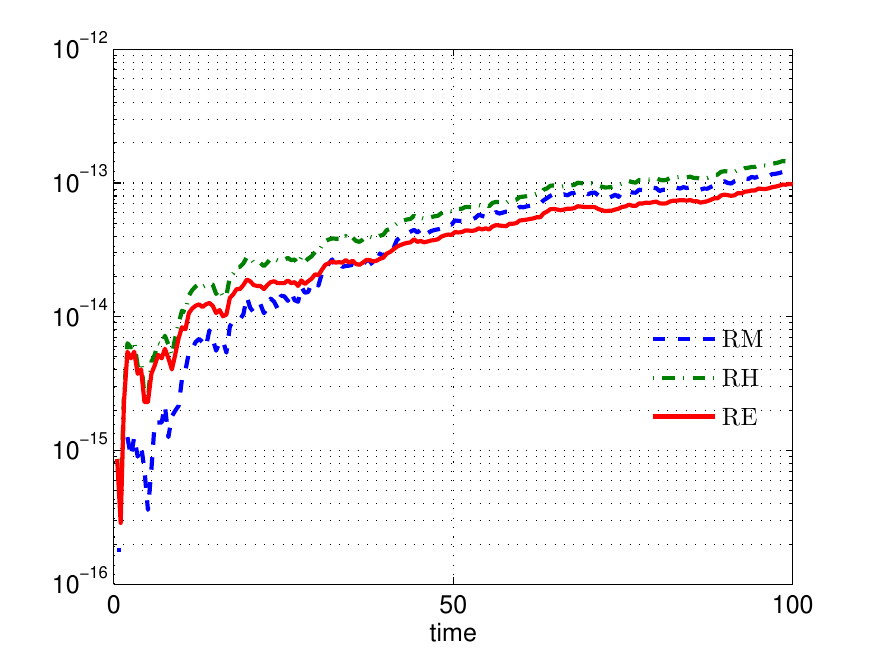}
  \caption*{(c)  4th-LI-PC-IF-RKM}
\end{minipage}
\caption{\footnotesize The relative residuals on three discrete conservation laws using the three proposed schemes for the 2D Schr\"odinger equation \eqref{NLS-equation}. We take the parameter  $\beta=-1$, time step $\tau=0.005$ and Fourier node $32\times 32$, and the initial condition is chosen as the exact solution at $t=0$.}\label{LI-EI-2d-scheme:fig:2}
\end{figure}

We choose the computational domain $\Omega=[0,2\pi]^2$ and take parameters $A=1,\ k_1=k_2=1$, $\beta=-1$ and $C_0=0$. We fix the Fourier node $32\times 32$ so that spatial errors can be neglected here. The $L^2$ errors and $L^{\infty}$ errors in numerical solution of $\phi$ at time $T=10$ are calculated by using the three  proposed schemes.  We take various time steps $\tau= \frac{5\times 10^{-3}}{2^i},\ i=0,1,2,3$ for 3rd-LI-IF-RKM and 4th-LI-IF-RKM, and for 4th-LI-PC-IF-RKM, we take various time steps $\tau=\frac{0.2}{2^i},\ i=0,1,2,3$. The results are summarized in Fig. \ref{LI-EI-2d-scheme:fig:1}. It is clear to see that 3rd-LI-IF-RKM is third-order accuracy in time, and 4th-LI-IF-RKM and  4th-LI-PC-IF-RKM have fourth-order accuracy in time, respectively. In Fig. \ref{LI-EI-2d-scheme:fig:2}, we investigate relative residuals  on the mass, Hamiltonian energy and quadratic energy using our schemes, which shows that the proposed schemes can preserve the discrete modified energy \eqref{LI-EI-F-energy} exactly, and the residuals on the discrete mass and Hamiltonian energy are very small.

Subsequently, we focus on computing the dynamics of turbulent superfluid governs by the following nonlinear Schr\"odinger equation
\begin{align} \label{NLS-equation2}
{\rm i} \partial_t \phi({\bf x},t) =-\frac{1}{2} \Delta \phi({\bf x},t) + \beta |\phi({\bf x},t)|^2 \phi({\bf x},t),  \  {\bf x} \in \mathbb R^d,\ d=2,3,\ t > 0,
\end{align}
with a initial data corresponding to a superfluid with a uniform condensate density and a phase which has a random spatial distribution~\cite{antoine15,Kobay}.

\begin{figure}[H]
    \centering
             \includegraphics[height=2.6cm,width=0.17 \linewidth]{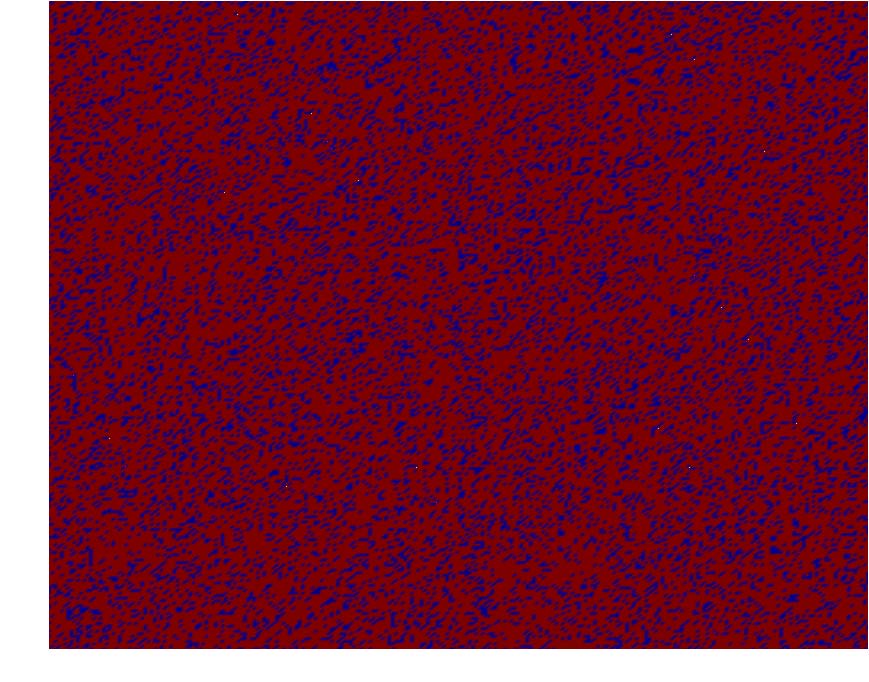}  \hspace{-3mm}
             \includegraphics[height=2.6cm,width=0.17 \linewidth]{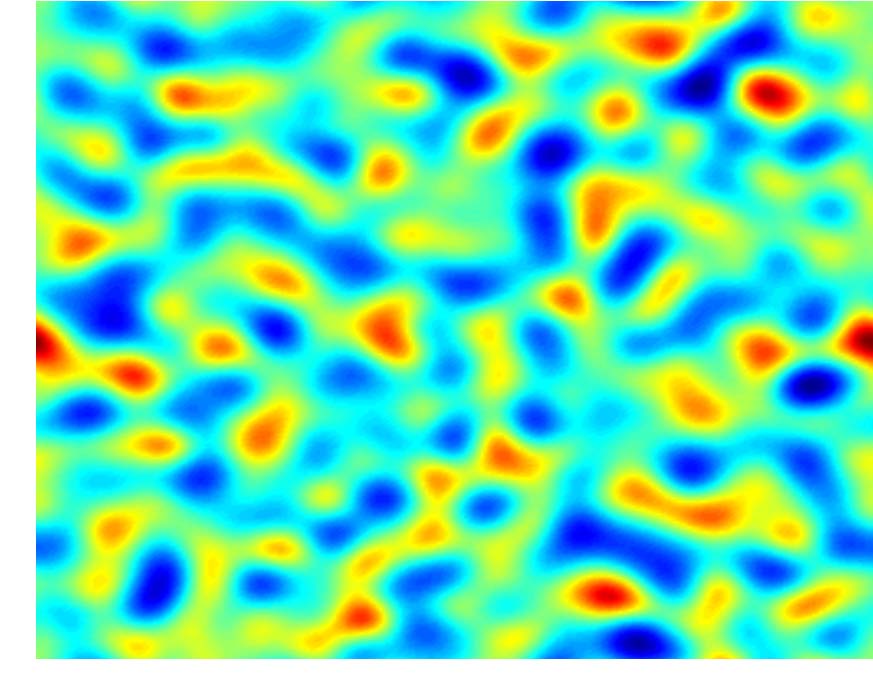} \hspace{-3mm}
	         \includegraphics[height=2.6cm,width=0.17 \linewidth]{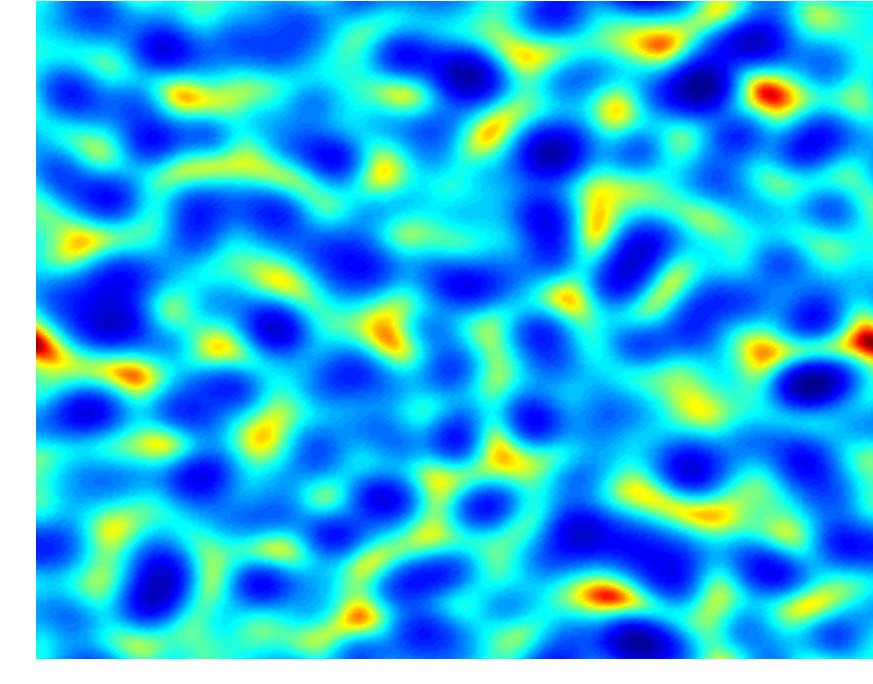} \hspace{-3mm}
             \includegraphics[height=2.6cm,width=0.17 \linewidth]{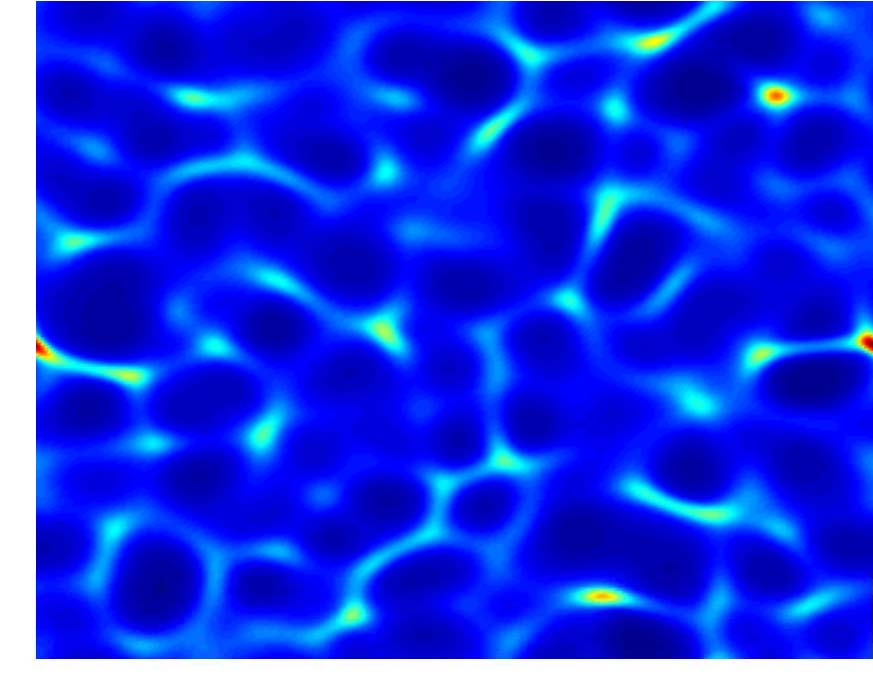} \hspace{-3mm}
	         \includegraphics[height=2.6cm,width=0.17 \linewidth]{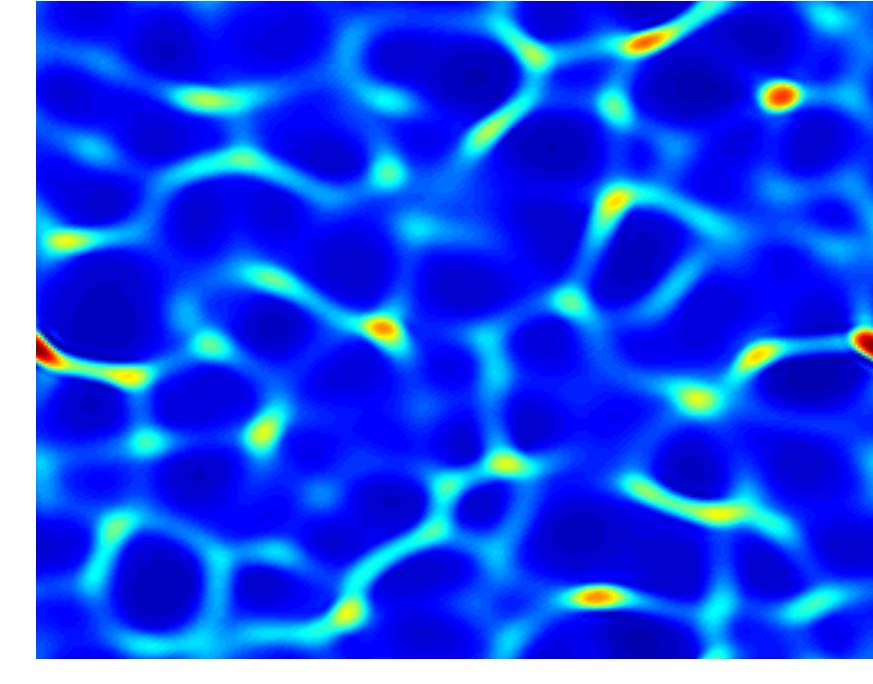}\\
             ~\includegraphics[height=2.6cm,width=0.17 \linewidth]{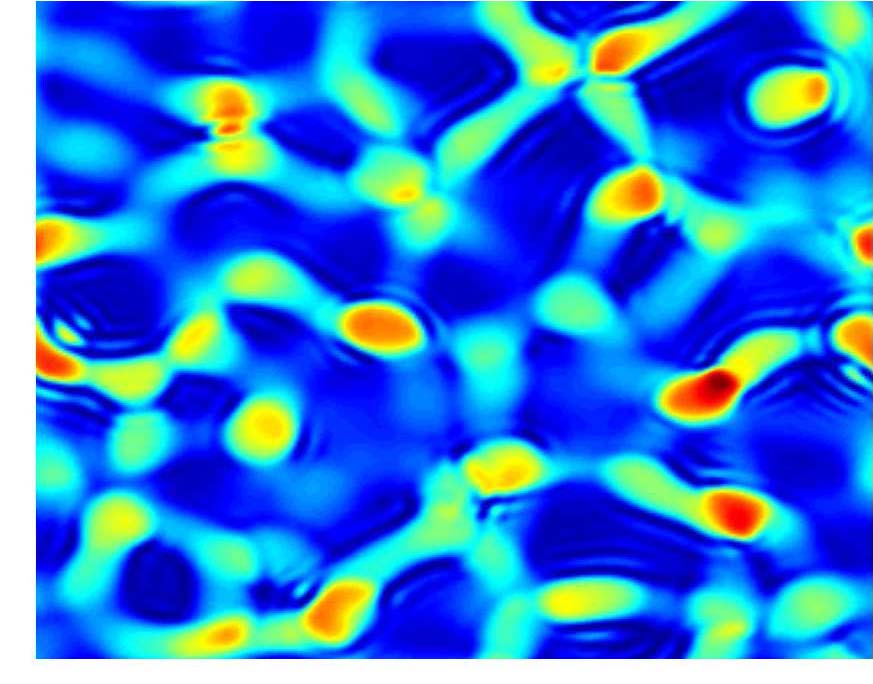} \hspace{-3mm}
	         \includegraphics[height=2.6cm,width=0.17 \linewidth]{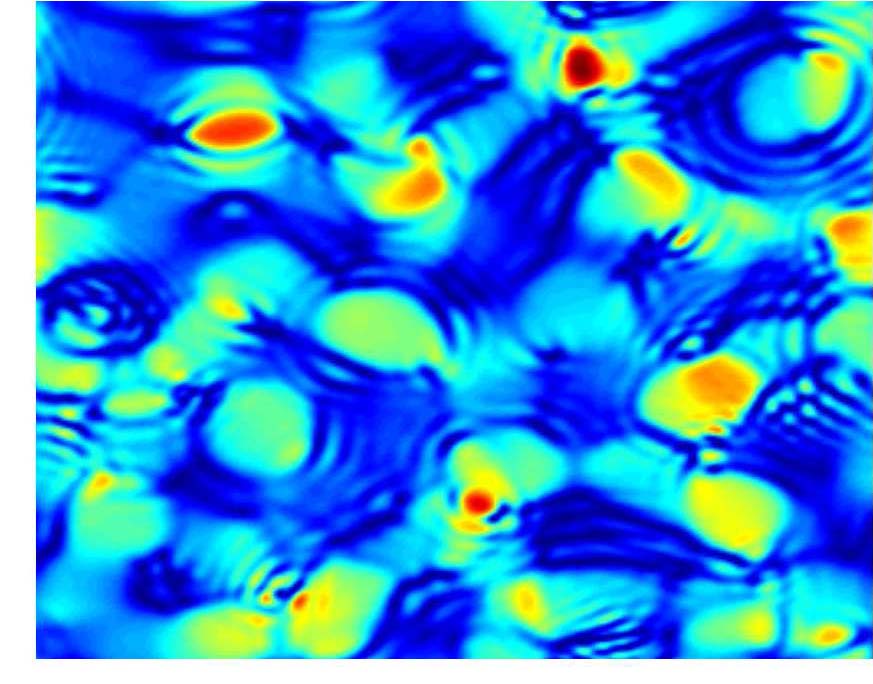} \hspace{-3mm}
             \includegraphics[height=2.6cm,width=0.17\linewidth]{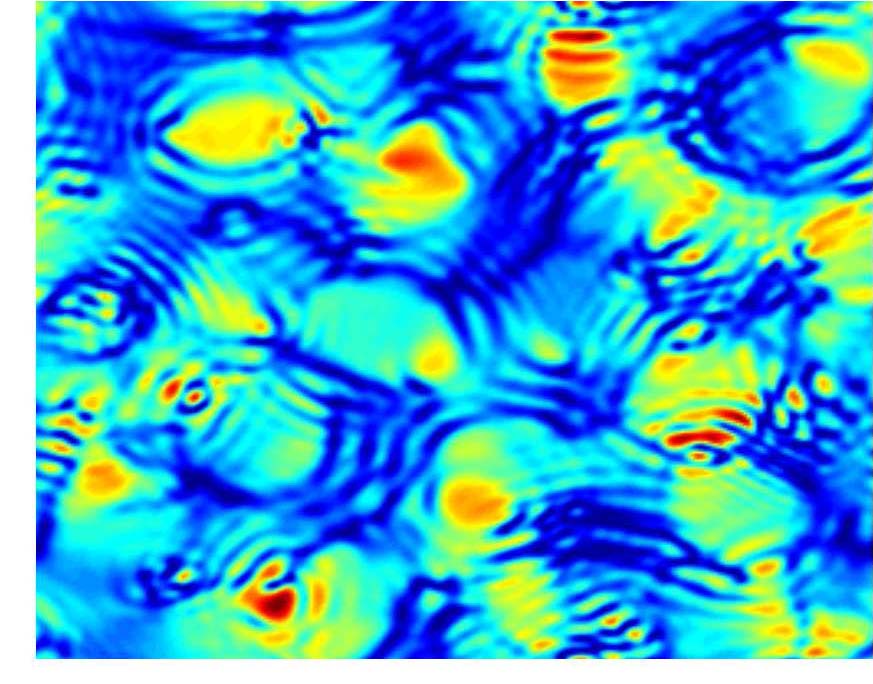} \hspace{-3mm}
             \includegraphics[height=2.6cm,width=0.17 \linewidth]{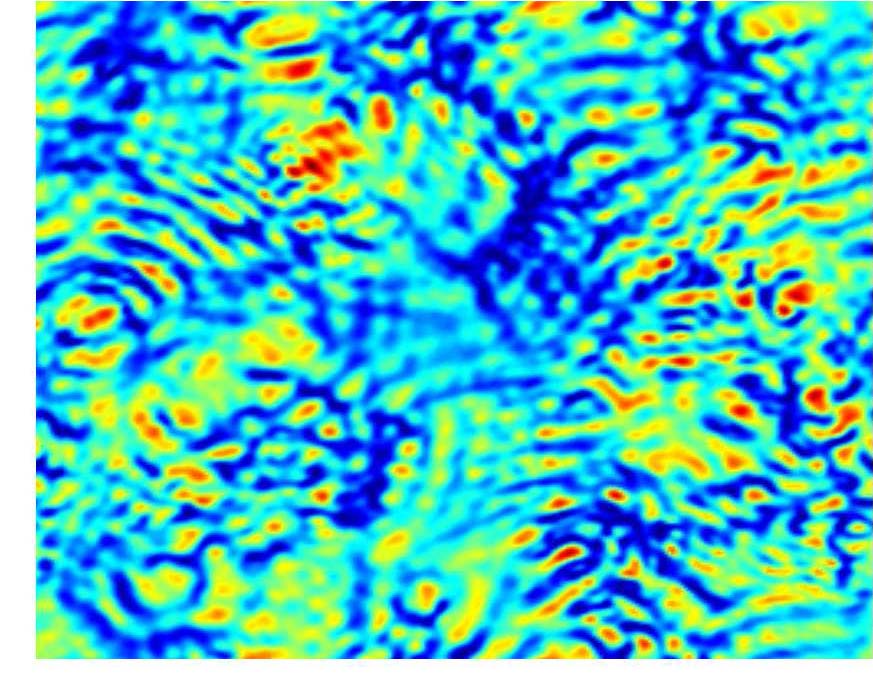} \hspace{-3mm}
             \includegraphics[height=2.6cm,width=0.17 \linewidth]{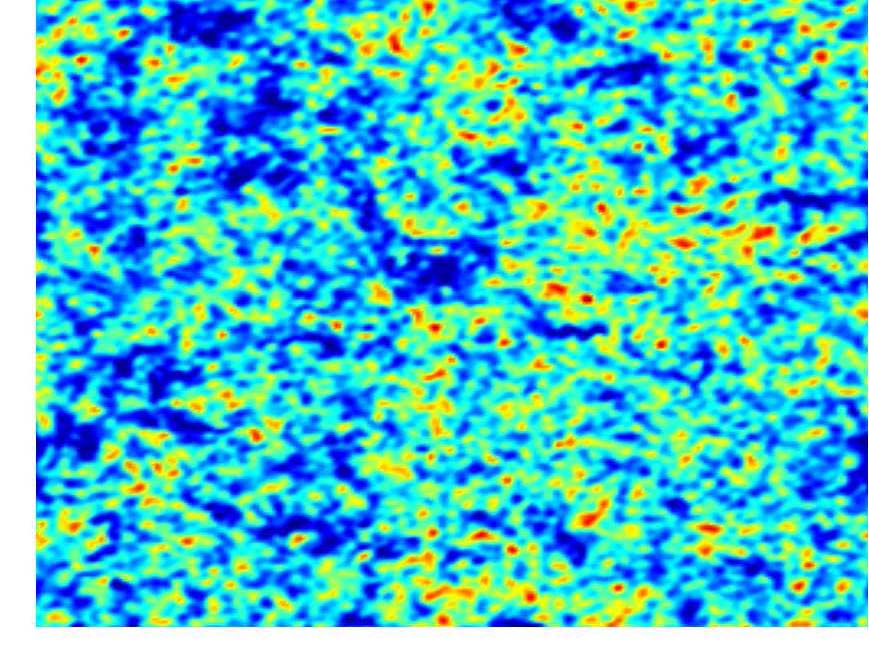}
\setlength{\abovecaptionskip}{-0.5mm}

\caption{\footnotesize { Modulus of the solution (i.e., $|\phi|^2$) for the dynamics of a turbulent
superfluid at different times $t=0,0.0005,0.002,0.006,0.01,0.2,0.3,0.4,0.7,5$ (in order from left to right and from top to bottom). We take the parameter $\beta=10$, time step $\tau=0.0001$ and Fourier node $256\times 256$, respectively.}}
\label{LI-EI-2d-scheme:fig:3}
\end{figure}

\begin{figure}[H]
\begin{minipage}[t]{45mm}
\includegraphics[width=45mm]{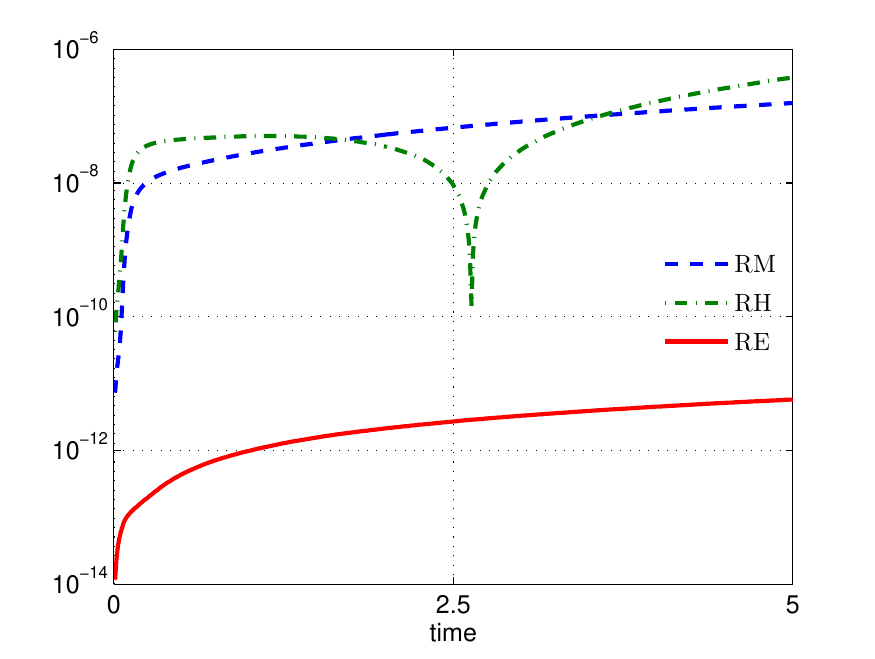}
  \caption*{(a)  3rd-LI-IF-RKM}
\end{minipage}
\begin{minipage}[t]{45mm}
\includegraphics[width=45mm]{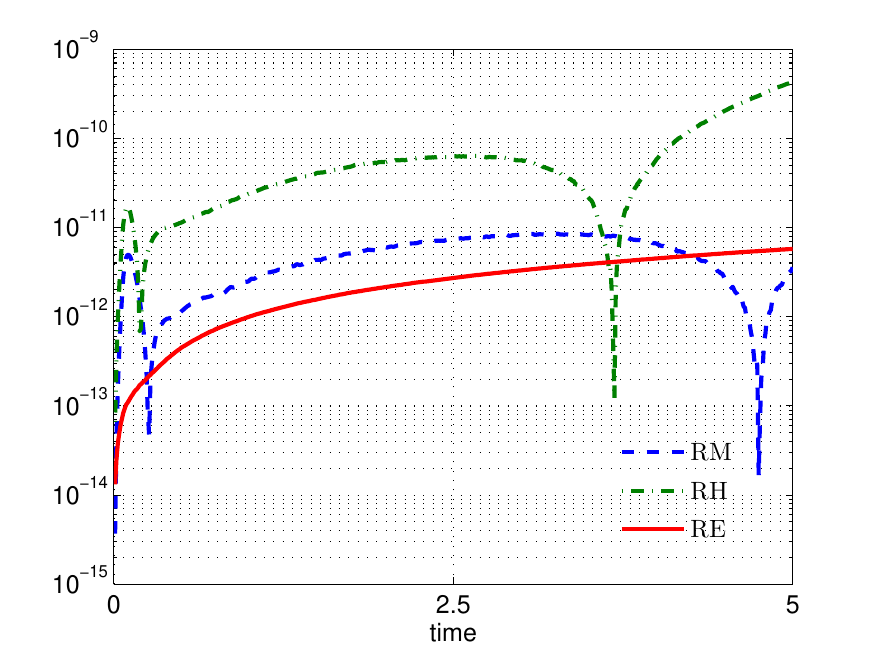}
  \caption*{(b)  4th-LI-IF-RKM}
\end{minipage}
\centering\begin{minipage}[t]{45mm}
\includegraphics[width=45mm]{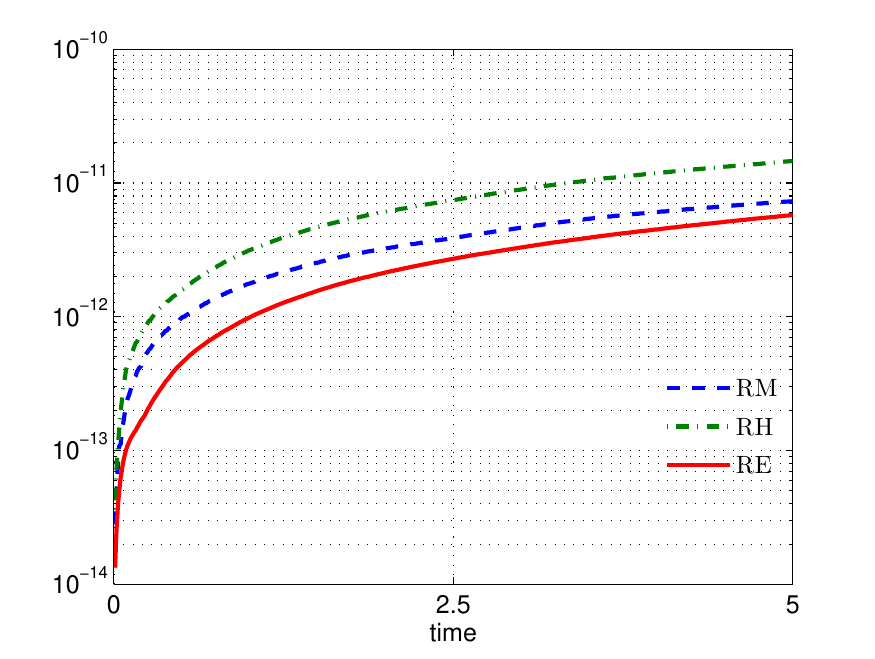}
  \caption*{(c)  4th-LI-PC-IF-RKM}
\end{minipage}
\caption{\footnotesize Relative  residuals on three discrete conservation laws using the three proposed numerical schemes for the 2D Schr\"odinger equation \eqref{NLS-equation2}. We  take the parameter $\beta=10$, time step $\tau=0.0001$ and Fourier node $256\times 256$, respectively.}
   \label{LI-EI-2d-scheme:fig:4}
\end{figure}

Following \cite{antoine15,Kobay}, we choose the computational domain $\Omega=[-10,10]^2$ with a periodic boundary condition. By using a similar procedure in~\cite{Kobay}, the initial data $\phi_0$ is chosen as $
\phi_0(x,y)=e^{{\rm i} \psi(x,y)}$, where $\psi$ is a random gaussian field with a covariance function given by $c(x,y)=e^{-(x^2+y^2)/2}$. {The modulus of the solution (i.e., $|\phi|^2$) for the dynamics of a turbulent
superfluid at different times provided by 3rd-LI-IF-RKM are summarized in Fig. \ref{LI-EI-2d-scheme:fig:3}.} We can clearly see that the phenomenons of wavy tremble are clearly observed at $t=0.3$, and the numerical results are in good agreement with those given in Refs.\cite{antoine15,Kobay}. Since the modulus of the solution computed by 4th-LI-IF-RKM and 4th-LI-PC-IF-RKM are similar to Fig. \ref{LI-EI-2d-scheme:fig:4}. For brevity, we omit it here.  In Fig. \ref{LI-EI-2d-scheme:fig:4} we show the relative residuals on the three conservation laws, which behaves similarly as that given in Fig. \ref{LI-EI-2d-scheme:fig:2}.

\vspace{-1mm}

\subsection{3D nonlinear Schr\"odinger equation}

In this subsection, we use the proposed schemes to solve the three dimensional (3D) nonlinear Schr\"odinger equation \eqref{NLS-equation}, which possesses the following analytical solution
\begin{align*}
\phi(x,y,z,t)=\text{exp}({\rm i}(k_1x+k_2y+k_3z-wt)), ~~w=k_1^2+k_2^2+k_3^2-\beta,\ (x,y,z)\in\Omega.
\end{align*}

\begin{figure}[H]
\centering\begin{minipage}[t]{60mm}
\includegraphics[width=60mm]{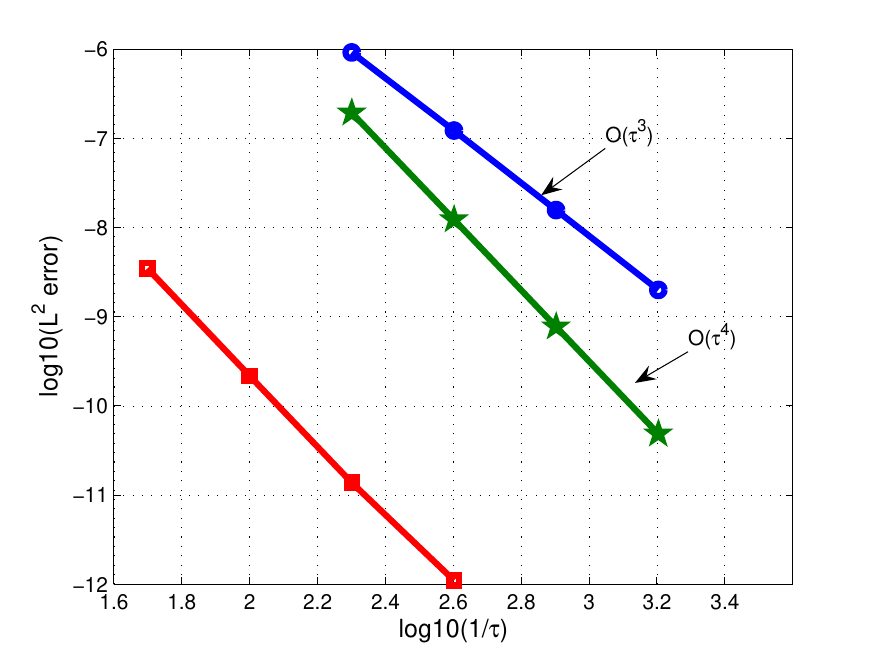}
\end{minipage}
\begin{minipage}[t]{60mm}
\includegraphics[width=60mm]{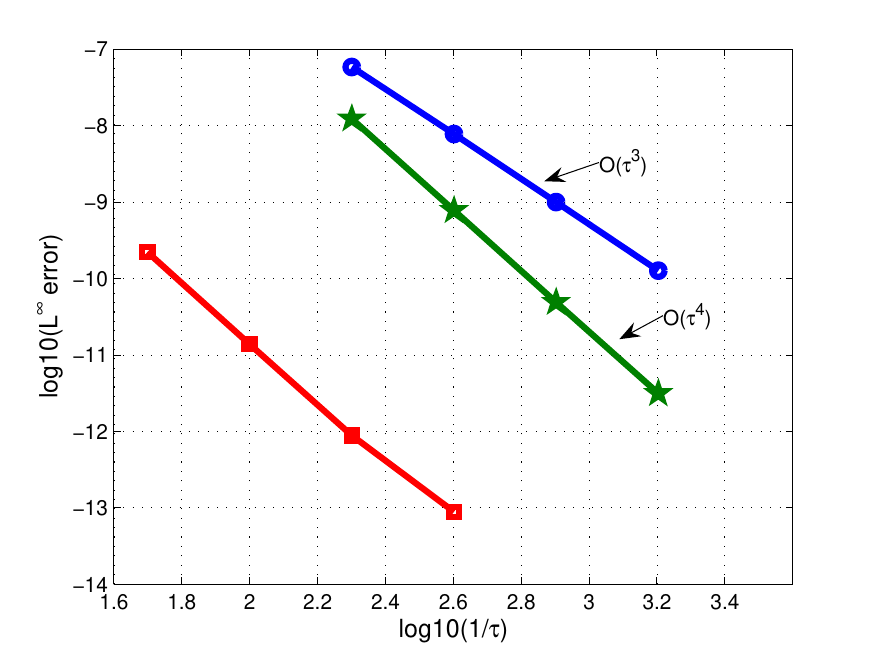}
\end{minipage}
\caption{ Time step refinement tests using the three proposed numerical methods for the 3D Schr\"odinger equation \eqref{NLS-equation}.}\label{LI-EI-3d-scheme:fig:1}
\end{figure}

\begin{figure}[H]
\begin{minipage}[t]{45mm}
\includegraphics[width=45mm]{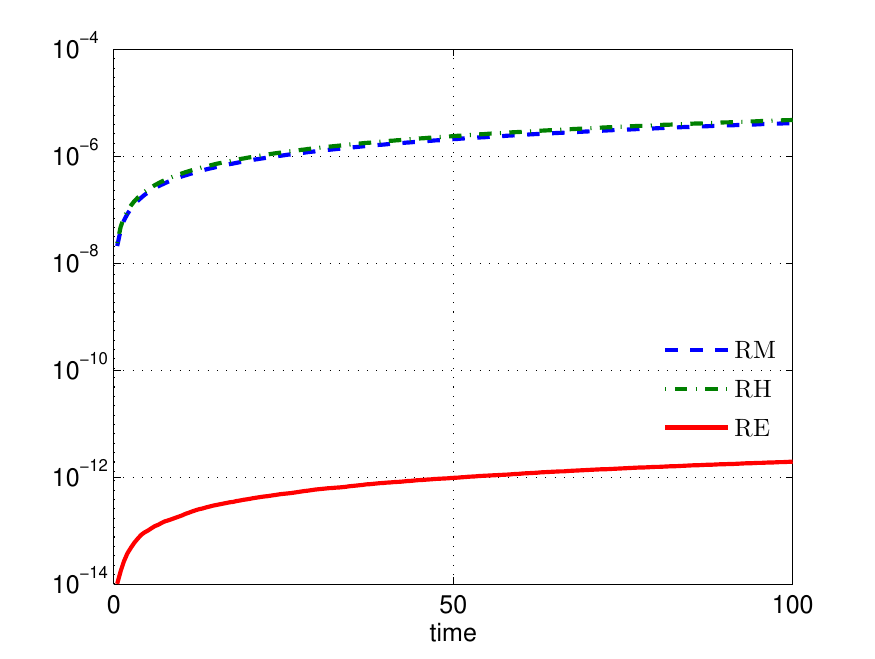}
  \caption*{(a)  3rd-LI-IF-RKM}
\end{minipage}
\begin{minipage}[t]{45mm}
\includegraphics[width=45mm]{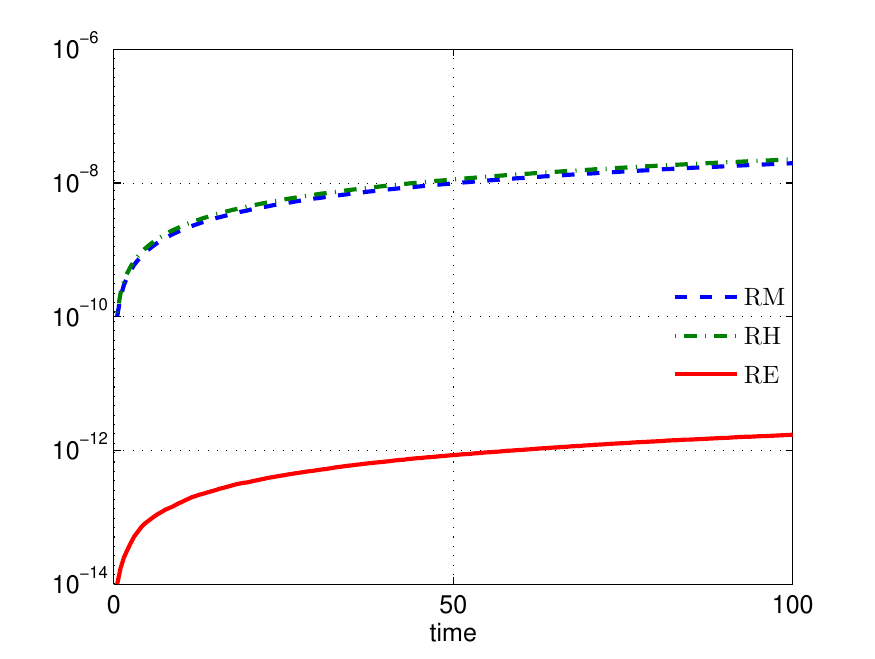}
  \caption*{(b)  4th-LI-IF-RKM}
\end{minipage}
\centering\begin{minipage}[t]{45mm}
\includegraphics[width=45mm]{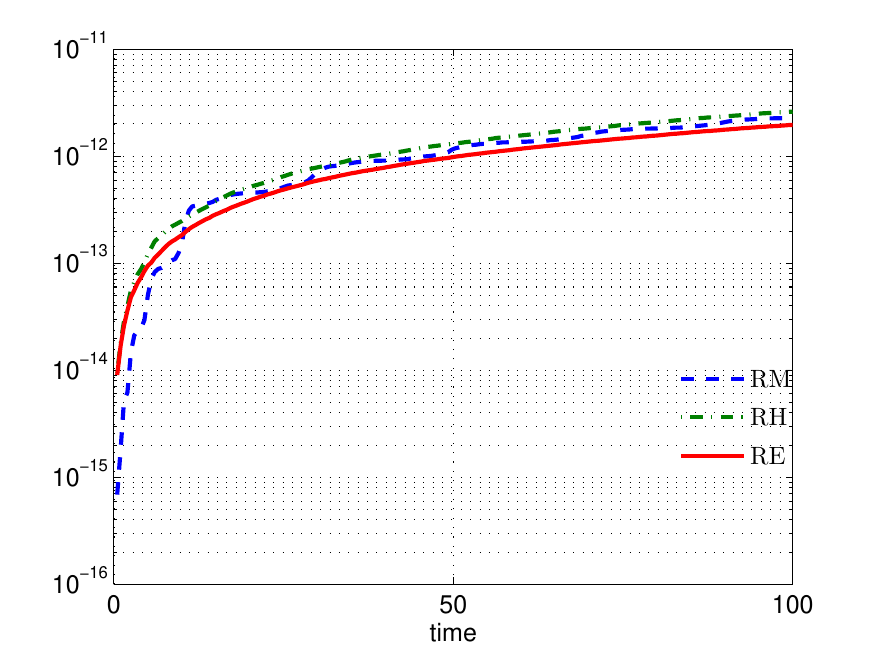}
  \caption*{(c)  4th-LI-PC-IF-RKM}
\end{minipage}
\caption{\footnotesize The relative residuals on three discrete conservation laws using the proposed numerical schemes for the 3D Schr\"odinger equation \eqref{NLS-equation}. We take the parameter $\beta=-1$, time step $\tau=0.005$ and Fourier node $16\times 16\times 16$, and the initial condition is chosen as the exact solution at $t=0$.}\label{LI-EI-3d-scheme:fig:2}
\end{figure}

We choose the computational domain $\Omega=[0,2\pi]^3$ and take parameters $\ k_1=k_2=k_3=1$, $\beta=-1$ and $C_0=0$. We fix the Fourier node $16\times 16\times 16$ so that spatial errors can be neglected here. The $L^2$ errors and $L^{\infty}$ errors in numerical solution of $\phi$ at time $T=1$ are calculated by using three proposed schemes. We take various time steps $\tau= \frac{5\times 10^{-3}}{2^i},\ i=0,1,2,3$ for 3rd-LI-IF-RKM and 4th-LI-IF-RKM, and for 4th-LI-PC-IF-RKM, we take various time steps $\tau=\frac{0.2}{2^i},\ i=0,1,2,3$. The results are summarized in Fig. \ref{LI-EI-3d-scheme:fig:1}. We can observe that 3rd-LI-IF-RKM has third-order accuracy in time, and  4th-LI-IF-RKM and 4th-LI-PC-IF-RKM have fourth-order accuracy in time, respectively. In Fig. \ref{LI-EI-3d-scheme:fig:2}, we investigate relative residuals on the mass, Hamiltonian energy
and quadratic energy by using 3rd-LI-IF-RKM, 4th-LI-IF-RKM and 4th-LI-PC-IF-RKM, respectively, on the time interval $t\in[0,100]$, which behaves similarly as those given in 1D and 2D cases.

\begin{figure}[H]
\centering
    \includegraphics[height=4.5cm,width=0.3\linewidth]{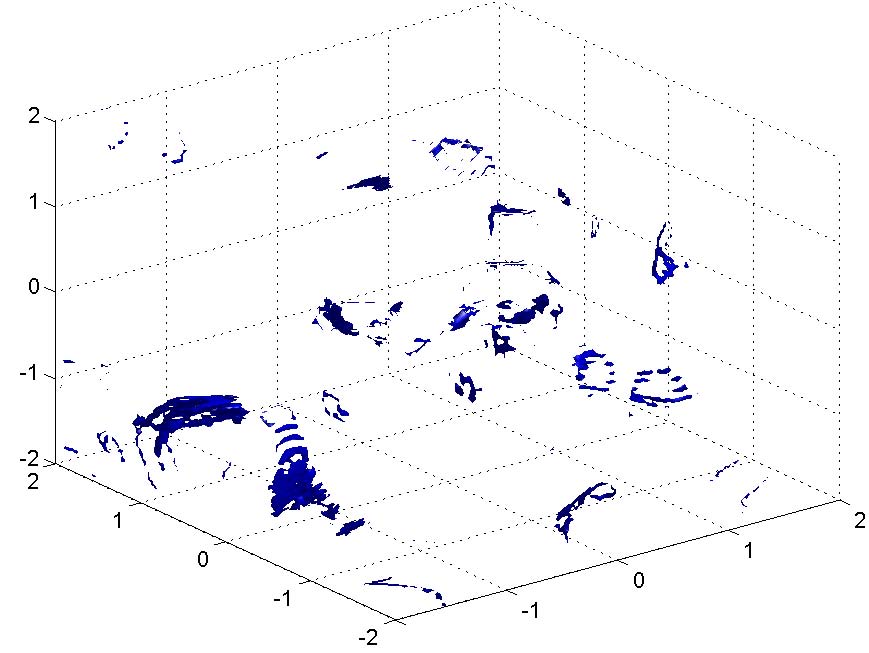}
    \includegraphics[height=4.5cm,width=0.3\linewidth]{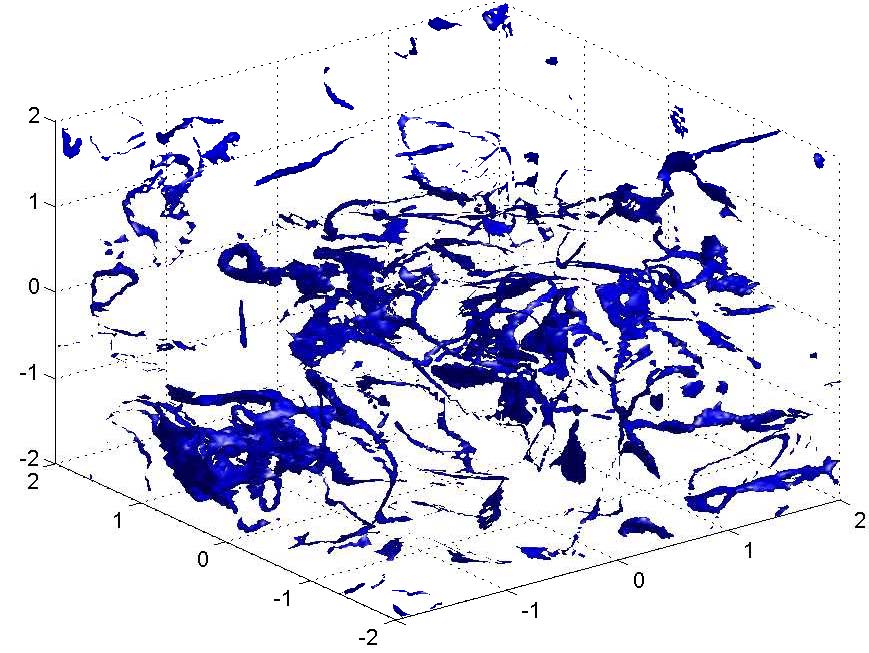}
	\includegraphics[height=4.5cm,width=0.3\linewidth]{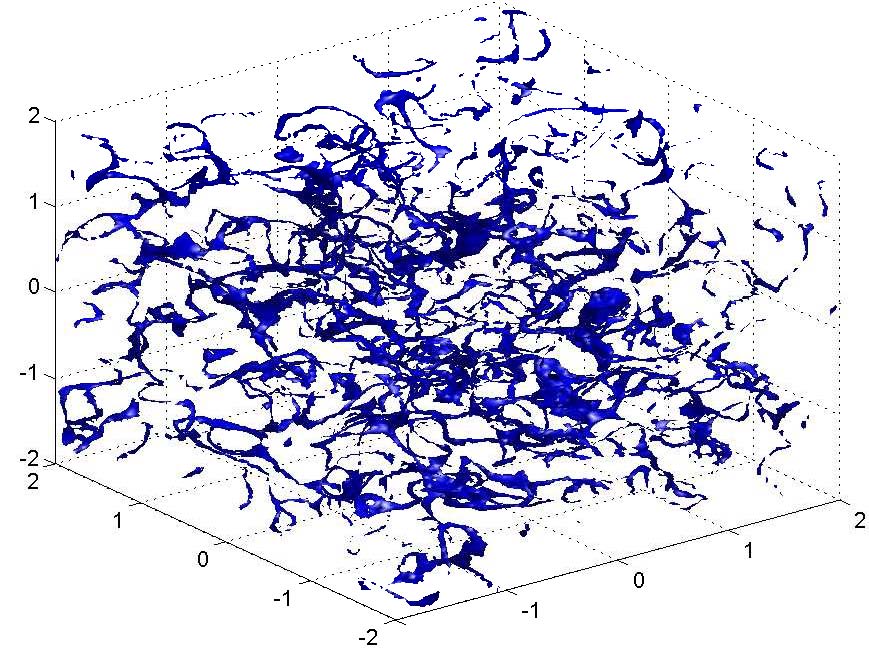}\\
	\includegraphics[height=4.5cm,width=0.3\linewidth]{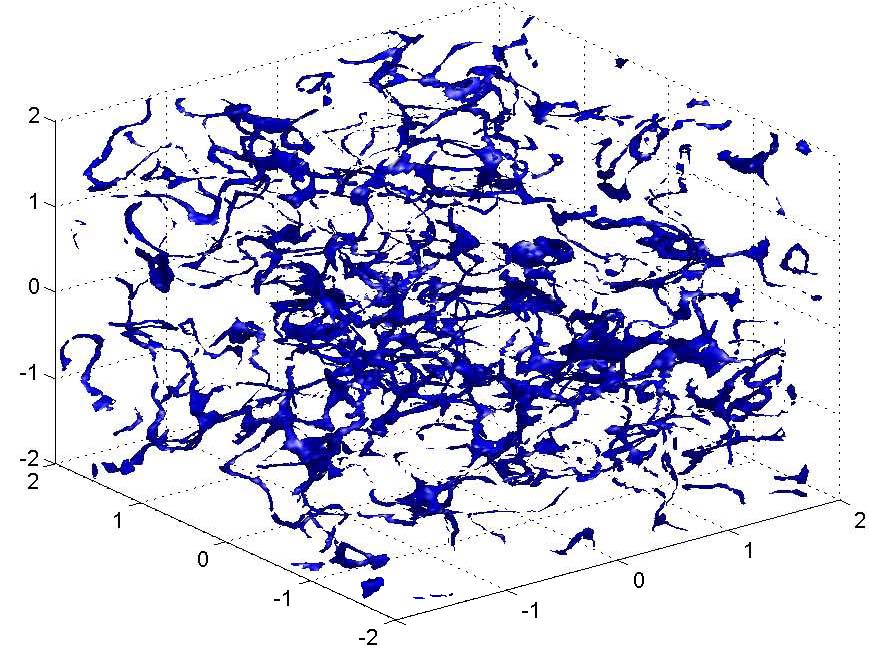}
    \includegraphics[height=4.5cm,width=0.3\linewidth]{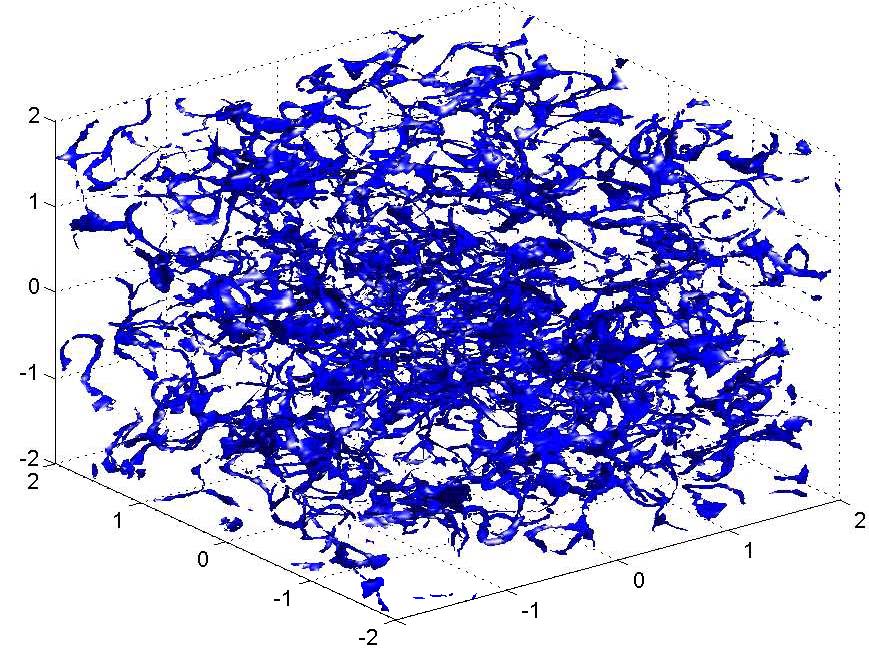}
	\includegraphics[height=4.5cm,width=0.3\linewidth]{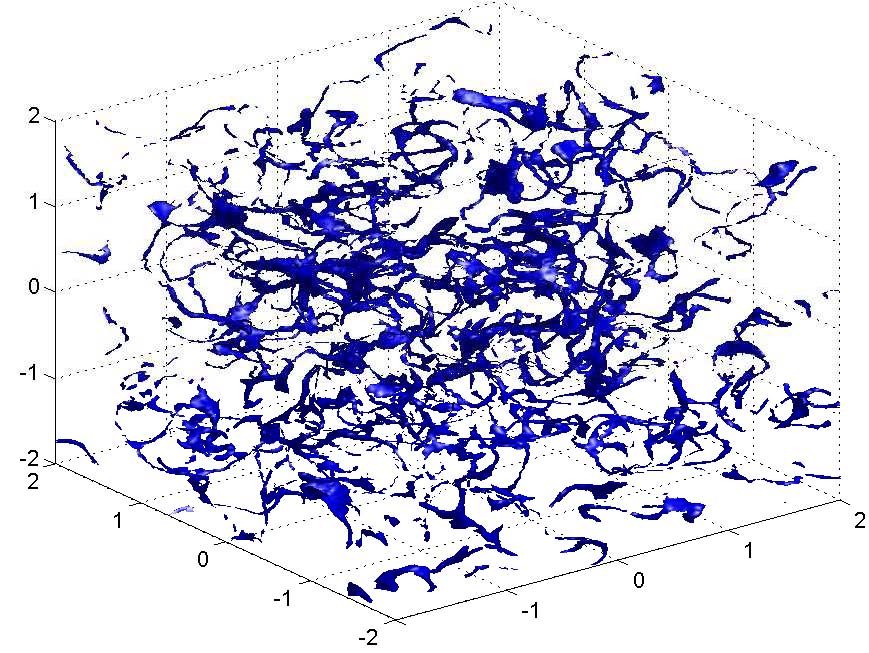}
	\caption{\small { Turbulence in a superfluid: modulus of the solution (i.e., $|\phi|^2$) at different times $t=0.04,0.06,0.5,0.6,0.8,1.0$ (in order from left to right and from top to bottom). We take the parameter $\beta=0.001$, time step $\tau=0.001$ and Fourier node $128\times 128\times 128$, respectively.}}
	\label{LI-EI-3d-scheme:fig:3}
\end{figure}

\begin{figure}[H]
\begin{minipage}[t]{45mm}
\includegraphics[width=45mm]{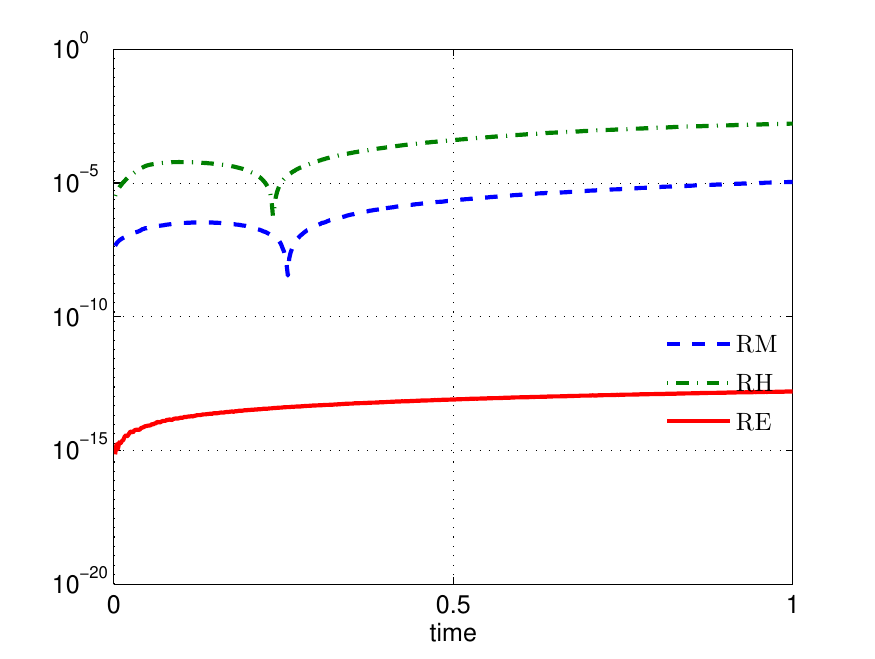}
  \caption*{(a)  3rd-LI-IF-RKM}
\end{minipage}
\begin{minipage}[t]{45mm}
\includegraphics[width=45mm]{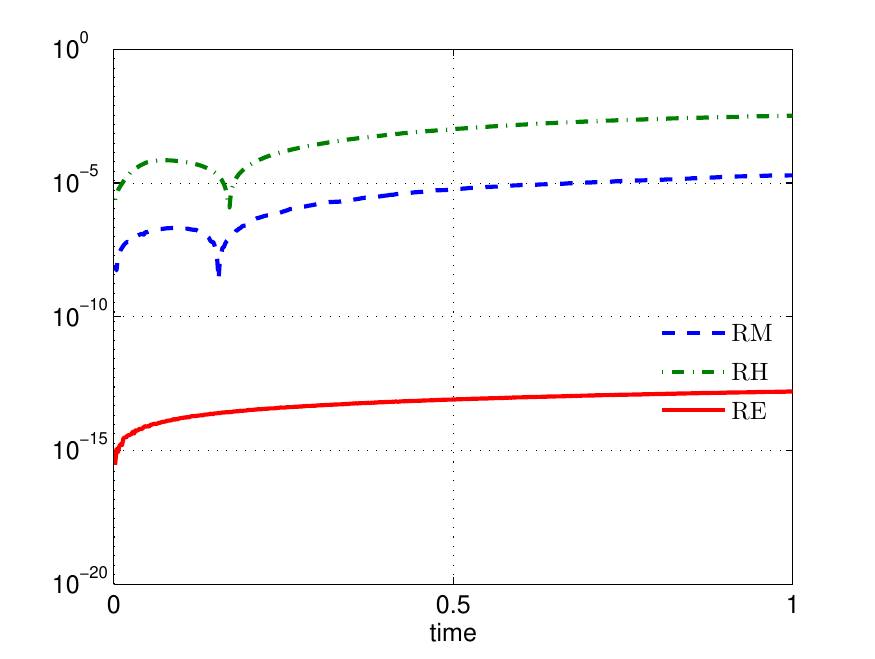}
  \caption*{(b)  4th-LI-IF-RKM}
\end{minipage}
\centering\begin{minipage}[t]{45mm}
\includegraphics[width=45mm]{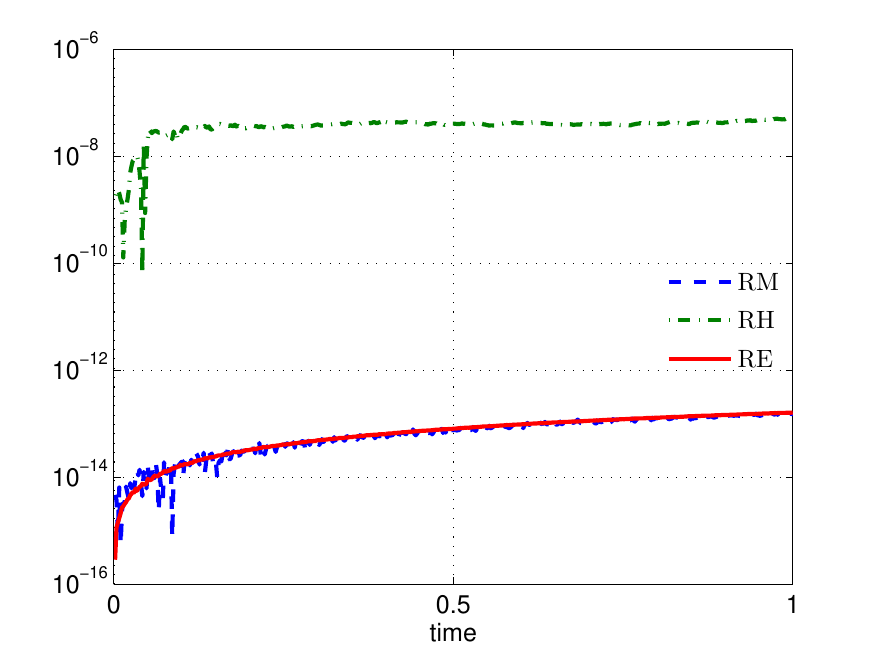}
  \caption*{(c)  4th-LI-PC-IF-RKM}
\end{minipage}
\caption{\footnotesize Relative residuals of three discrete conservation laws using the three proposed numerical schemes for the Schr\"odinger equation \eqref{NLS-equation2}. We take the parameter $\beta=0.001$, time step $\tau=0.001$ and Fourier node $128\times 128\times 128$, respectively. }\label{LI-EI-3d-scheme:fig:4}
\end{figure}

Finally, we use our schemes to investigate the dynamics of the 3D nonlinear Schr\"odinger equation~\eqref{NLS-equation2} (i.e., $d=3$). Inspired by \cite{antoine15,Kobay},   the initial data is chosen to be  associated with a superfluid with a uniform density and a random phase, and we take the computational domain $\Omega=[-2,2]^3$ with a periodic boundary condition. { The modulus of the solution (i.e., $|\phi|^2$) for the formation of a 3D turbulent phenomenon at the different times provided by 4th-LI-IF-RKM are summarized in Fig.~\ref{LI-EI-3d-scheme:fig:3}}. It can be observed that the vortices filamentation forms gradually during the time evolution. In particular, at time $t=0.8$, the massive vortices filamentation filled to the whole spatial domain, which is consistent with the numerical results presented in~\cite{antoine15}. Here, we omit the numerical results computed by 3rd-LI-IF-RKM and 4th-LI-PC-IF-RKM, since they are similar to Fig. \ref{LI-EI-3d-scheme:fig:3}. The relative residuals on the three conservation laws provided by the proposed numerical schemes are displayed in Fig. \ref{LI-EI-3d-scheme:fig:4}. As is illustrated in the figure, the proposed schemes can preserve the modified energy exactly and conserve the discrete mass and Hamiltonian energy well.

\section{Conclusions}\label{Sec:LI-EP:5}

In this paper, based on the Lawson transformation, the Runge-Kutta method and the idea of the SAV approach, we develop a novel class of
high-order linearly implicit energy-preserving integrating
factor Runge-Kutta methods for the nonlinear Schr\"odinger equation. The construction
of the proposed schemes are two key steps, as follows: (i) We first obtain a linearized system which satisfies a quadratic energy by using the extrapolation technique
or prediction-correction strategy to the nonlinear terms of the SAV reformulation; (ii) Then, based on the integrating factor Runge-Kutta method under certain circumstances for the RK coefficients (see \eqref{LI-EI-coefficient}), a class of linearly implicit energy-preserving schemes are obtained. The proposed scheme can reach high-order accuracy in time and inherit the advantage of the ones provided by the classical SAV approach. Extensive numerical examples are addressed to illustrate
the efficiency and accuracy of our new schemes. Compared with some existing structure-preserving schemes, the proposed scheme based on the prediction-correction strategy shows remarkable superiorities. {However, we should note that the proposed scheme can only preserve a modified energy and we are sure that the high-order linearly implicit integrating
factor Runge-Kutta methods which can preserve original energy will be more advantages than the proposed schemes in the robustness and the long-time accurate computations. Thus, an interesting topic for future studies is whether it is possible to construct high-order linearly implicit schemes which
can preserve the original energy.}

\section*{Acknowledgments}
The authors would like to express sincere gratitude to the referees for their insightful comments and suggestions.
This work is supported by the National Natural Science Foundation of China (Grant Nos. 11901513, 11971481, 12071481), the National Key R\&D Program of China (Grant No. 2020YFA0709800), the Natural Science Foundation of Hunan (Grant Nos. 2021JJ40655, 2021JJ20053), the Yunnan Fundamental Research Projects (Nos. 202101AT070208, 202001AT070066,
202101AS070044), the High Level Talents Research Foundation Project of Nanjing Vocational College of Information Technology (Grant No. YB20200906), and the Foundation of Jiangsu Key Laboratory for Numerical Simulation of Large Scale Complex Systems (Grant No. 202102).
 Jiang and Cui are in particular grateful to Prof. Yushun Wang and Dr. Yuezheng Gong for fruitful discussions.


\end{document}